\documentclass[a4paper]{article}
\usepackage{amsthm,amsmath,amssymb}
\newtheorem{teor}{Theorem}[section]
\newtheorem{defin}[teor]{Definition}
\newtheorem{lemm}[teor]{Lemma}
\newtheorem{osse}[teor]{Remark}
\newtheorem{prop}[teor]{Proposition}
\newtheorem{defi}[teor]{Definition}
\newtheorem{coro}[teor]{Corollary}
\newtheorem{prob}[teor]{Problem}

\newcommand{\bele}{\begin{lemm}\begin{sl}}
\newcommand{\enle}{\end{sl}\end{lemm}}
\newcommand{\bedef}{\begin{defi}\begin{sl}}
\newcommand{\eddef}{\end{sl}\end{defi}}
\newcommand{\bete}{\begin{teor}\begin{sl}}
\newcommand{\ente}{\end{sl}\end{teor}}
\newcommand{\beos}{\begin{osse}\begin{rm}}
\newcommand{\eddos}{\end{rm}\end{osse}}
\newcommand{\bepr}{\begin{prop}\begin{sl}}
\newcommand{\empr}{\end{sl}\end{prop}}
\newcommand{\bepro}{\begin{prob}\begin{rm}}
\newcommand{\empro}{\end{rm}\end{prob}}
\newcommand{\bede}{\begin{defin}\begin{sl}}
\newcommand{\edde}{\end{sl}\end{defin}}
\newcommand{\beco}{\begin{coro}\begin{sl}}
\newcommand{\enco}{\end{sl}\end{coro}}

\newcommand{\quext}{\quad\text}
\newcommand{\qquext}{\qquad\text}
\newcommand{\de}{\partial}

\newcommand{\RR}{\mathbb{R}}

\newcommand{\beeq}[1]{\begin{equation}\label{#1}}
\newcommand{\eddeq}{\end{equation}}

\newcommand{\beeqa}[1]{\begin{eqnarray}\label{#1}}
\newcommand{\eddeqa}{\end{eqnarray}}

\newcommand{\beal}[1]{\begin{align}\label{#1}}
\newcommand{\eddal}{\end{align}}

\newcommand{\bespl}[1]{\begin{split}\label{#1}}
\newcommand{\edspl}{\end{split}}

\newcommand{\bega}[1]{\begin{gather}\label{#1}}
\newcommand{\edga}{\end{gather}}

\newcommand{\beeqax}{\begin{eqnarray*}}
\newcommand{\eddeqax}{\end{eqnarray*}}

\def\qed{\ifmmode 
  \else \leavevmode\unskip\penalty9999 \hbox{}\nobreak\hfill
  \fi
  \quad\hbox{\hskip.5em\vrule width.4em height.6em depth.05em\hskip.1em}}
\def\endproofsym{\qed}
\renewenvironment{proof}[1][Proof]{\trivlist\item[\hskip\labelsep{\hskip0pt
    {\normalfont\scshape#1.}\hskip .321429\parindent}]\ignorespaces}
{\endproofsym\endtrivlist}
\def\endnobox{\def\endproofsym{}\end{proof}\def\endproofsym{\qed}}

\newcommand{\no}{\nonumber}

\newcommand{\beeqao}{\begin{eqnarray}\no}
\newcommand{\bealo}{\begin{align}\no}
\newcommand{\besplo}{\begin{split}\no}
\newcommand{\begao}{\begin{gather}\no}

\newcommand{\implica}{\quad\Longrightarrow\quad}

\newcommand{\duav}[1]{\langle{#1}\rangle}

\newcommand{\duavg}[1]{\left\langle{#1}\right\rangle}
\newcommand{\duavb}[1]{\big\langle{#1}\big\rangle}

\newcommand{\perogni}{\forall\,}
\newcommand{\esiste}{\exists\,}
\newcommand{\itt}{\int_0^t}

\newcommand{\io}{\int_\Omega}
\newcommand{\ibaro}{\int_{\overline{\Omega}}}

\newcommand{\iga}{\int_\Gamma}

\newcommand{\epsi}{\varepsilon}

\newcommand{\dd}{_\delta}
\newcommand{\ddm}{_{\delta,m}}
\newcommand{\dda}{_{\delta,a}}
\newcommand{\dds}{_{\delta,s}}

\newcommand{\OO}{_{\Omega}}

\newcommand{\fzw}{f_{0,w}}

\newcommand{\bn}{\boldsymbol{n}}

\newcommand{\dn}{\partial_{\bn}}

\newcommand{\lhs}{left hand side}
\newcommand{\rhs}{right hand side}

\DeclareMathOperator{\deriv}{d}

\DeclareMathOperator{\Id}{Id}

\DeclareMathOperator{\const}{const}

\newcommand{\HUH}{H^1(0,T;H)}

\newcommand{\HUVp}{H^1(0,T;V')}
\newcommand{\CZH}{C^0([0,T];H)}

\newcommand{\LDH}{L^2(0,T;H)}
\newcommand{\LDV}{L^2(0,T;V)}
\newcommand{\LDVp}{L^2(0,T;V')}

\newcommand{\LIV}{L^\infty(0,T;V)}
\newcommand{\LIW}{L^\infty(0,T;W)}

\newcommand{\LDW}{L^2(0,T;W)}

\let\TeXchi\chi
\def\chi{{\setbox0 \hbox{\mathsurround0pt
$\TeXchi$}\hbox{\raise\dp0 \copy0 }}}

\newcommand{\calX}{{\mathcal X}}

\newcommand{\calT}{{\mathcal T}}

\newcommand{\calU}{{\mathcal U}}
\newcommand{\calA}{{\mathcal A}}
\newcommand{\calK}{{\mathcal K}}
\newcommand{\calF}{{\mathcal F}}
\newcommand{\calE}{{\mathcal E}}
\newcommand{\calJ}{{\mathcal J}}

\newcommand{\calM}{{\mathcal M}}

\newcommand{\calB}{{\mathcal B}}

\newcommand{\calS}{{\mathcal S}}

\newcommand{\calW}{{\mathcal W}}

\newcommand{\zzd}{_{0,\delta}}
\newcommand{\zzs}{_{0,\sigma}}
\newcommand{\ssi}{_\sigma}
\newcommand{\ssid}{_{\sigma,\delta}}

\newcommand{\barO}{\overline{\Omega}}

\newcommand{\baru}{\overline{u}}
\newcommand{\barw}{\overline{w}}
\newcommand{\barf}{\overline{f}}

\newcommand{\agiu}{\underline{a}}
\newcommand{\asu}{\overline{a}}

\newcommand{\dit}{\deriv\!t}
\newcommand{\dis}{\deriv\!s}
\newcommand{\dix}{\deriv\!x}

\newcommand{\dixi}{\deriv\!\xi}
\newcommand{\dixis}{\deriv\!|\xi\dds|}
\newcommand{\dixist}{\deriv\!|\xi\dds(t)|}
\newcommand{\diT}{\deriv\!T}
\newcommand{\dir}{\deriv\!r}

\newcommand{\ddt}{\frac{\deriv\!{}}{\dit}}

\newcommand{\calL}{{\cal L}}

\numberwithin{equation}{section}
\begin{document}

\title{On a class of Cahn-Hilliard models
  with nonlinear diffusion}

\author{Giulio Schimperna\\
Dipartimento di Matematica, Universit\`a di Pavia,\\
Via Ferrata~1, I-27100 Pavia, Italy\\
E-mail: {\tt giusch04@unipv.it}\\
\and
Irena Paw\l ow\\
Systems Research Institute,\\
Polish Academy of Sciences\\
\mbox{}~~and Institute of Mathematics and Cryptology,\\
Cybernetics Faculty,\\
Military University of Technology,\\
S.~Kaliskiego 2, 00-908 Warsaw, Poland\\ 
E-mail: {\tt Irena.Pawlow@ibspan.waw.pl}
}


\maketitle
\begin{abstract}
 In the present work, we address a class of Cahn-Hilliard equations
 characterized by a nonlinear diffusive dynamics
 and possibly containing an additional sixth order term.
 This model describes the separation properties
 of oil-water mixtures, when a substance enforcing the mixing of the
 phases (a surfactant) is added. However, the model is also closely 
 connected with other Cahn-Hilliard-like equations
 relevant in different types of applications.
 We first discuss the existence of a weak solution to the sixth-order model 
 in the case when the configuration potential of the system has 
 a singular (e.g., logarithmic) character. Then, we study the behavior
 of the solutions in the case when the sixth order term 
 is let tend to 0, proving convergence to solutions of the 
 fourth order system in a special case. The fourth order system
 is then investigated by a direct approach and existence of
 a weak solution is shown under very general conditions 
 by means of a fixed point argument.
 Finally, additional properties of the solutions, like uniqueness
 and parabolic regularization, are discussed, both for the sixth order
 and for the fourth order model, under more restrictive assumptions on 
 the nonlinear diffusion term.
\end{abstract}


\vspace{2mm}

\noindent {\bf AMS (MOS) subject clas\-si\-fi\-ca\-tion:}%
~~35K35, 35K55, 35A01, 47H05.

\vspace{2mm}


\section{Introduction}
\label{sec:intro}

This paper is devoted to the mathematical analysis of the following
class of parabolic systems:
\begin{align}\label{CH1}
  & u_t - \Delta w = 0,\\
 \label{CH2}
  & w = \delta \Delta^2 u - a(u) \Delta u - \frac{a'(u)}2 |\nabla u|^2 
   + f(u) + \epsi u_t,
\end{align}
on $(0,T) \times \Omega$, $\Omega$ being a bounded 
smooth subset of $\RR^3$ and $T>0$ an assigned final
time. The restriction to the three-dimensional setting is
motivated by physical applications. Similar, or better results,
are expected to hold in space dimensions 1 and 2.
The system is coupled with
the initial and boundary conditions
\begin{align}\label{iniz-intro}
  & u|_{t=0} = u_0,
  \quext{in }\,\Omega,\\
 \label{neum-intro}
  & \dn u = \dn w = \delta \dn \Delta u = 0,
  \quext{on }\,\partial\Omega,\
  \quext{for }\,t\in(0,T)
\end{align}
and represents a variant of the 
Cahn-Hilliard model for phase separation in binary materials. 
The function $f$ stands for the derivative 
of a {\sl singular}\/ potential $F$ of a {\sl double obstacle}\
type. Namely, $F$ is assumed to be $+\infty$ outside a bounded interval
(assumed equal to $[-1,1]$ for simplicity), where the extrema
correspond to the pure states. A physically significant 
example is given by the so-called 
Flory-Huggins logarithmic potential  
\begin{equation}\label{logpot}
  F(r)=(1-r)\log(1-r)+(1+r)\log(1+r) - \frac\lambda2 r^2,
   \quad \lambda\ge 0.
\end{equation}
As in this example, we will assume $F$ to be at least
{\sl $\lambda$-convex}, i.e., convex up to a quadratic
perturbation. In this way, we can also allow for singular
potentials having more than two minima in the interval $[-1,1]$
(as it happens in the case of the oil-water-surfactant models
described below, where the third minimum appears in relation
to the so-called ``microemulsion'' phase).

We assume the coefficients $\delta,\epsi$ to be $\ge 0$, with the 
case $\delta>0$ giving rise to a {\sl sixth order}\ model
and the case $\epsi>0$ related to possible viscosity effects
that are likely to appear in several models of Cahn-Hilliard type
(see, e.g., \cite{Gu}). 
The investigation of the limits as $\delta$ or $\epsi$ tend to zero
provides validation of these models as the approximates of the limit fourth order
model.

The main novelty of system \eqref{CH1}-\eqref{CH2}
is related to the presence of the 
nonlinear function $a$ in \eqref{CH2}, which is supposed
smooth, bounded, and strongly positive (i.e., 
everywhere larger than some constant $\agiu>0$). Mathematically,
the latter is an unavoidable assumption as we are mainly interested
in the behavior of the problem when $\delta$ is let 
tend to $0$ and in the properties of the (fourth order)
limit system $\delta=0$. On the other hand, at least in the 
physical context of the sixth order model, it would also be
meaningful to admit $a$ to take negative values, as it may
happen in presence of the ``microemulsion'' phase 
(see \cite{GK93a,GK93b}). We will not deal with this 
situation, but we just point out that, as long as $\delta>0$
is fixed, this should create no additional mathematical 
difficulties since the nonlinear diffusion
term is then dominated by the sixth order term.

From the analytical point of view,
as a basic observation we can notice that this class of systems
has an evident variational structure. Indeed, (formally) 
testing \eqref{CH1} by $w$, \eqref{CH2} by $u_t$,
taking the difference of the obtained relations, 
integrating with respect to~space variables,
using the {\sl no-flux}\/ conditions \eqref{neum-intro},
and performing suitable integrations
by parts, one readily gets the {\sl a priori}\/ bound
\begin{equation}\label{energyineq}
  \ddt\calE\dd(u) + \| \nabla w \|_{L^2(\Omega)}^2 
   + \epsi \| u_t \|_{L^2(\Omega)}^2 
   = 0,
\end{equation}
which has the form of an {\sl energy equality} for the 
{\sl energy functional}
\begin{equation}\label{defiE}
  \calE\dd(u)=\io \Big( \frac\delta2 |\Delta u|^2
   + \frac{a(u)}2 |\nabla u|^2 
   + F(u) \Big),
\end{equation}
where the interface (gradient) part contains the nonlinear 
function $a$. In other words, the system \eqref{CH1}-\eqref{CH2}
arises as the $(H^1)'$-gradient flow problem for the 
functional $\calE\dd$.
While the literature on the fourth order Cahn-Hilliard model
with logarithmic free energy is very wide (starting from the 
pioneering work \cite{DD} up to more recent works like,
e.g., \cite{AW,GMS,MZ}, see also the recent review \cite {ChMZ} and
the references therein),
it seems that potentials of logarithmic types have 
never been considered in the case of a nonconstant 
coefficient $a$.
Similarly, the sixth order Cahn-Hilliard type equations, which appear as
models of various physical phenomena and 
have recently attracted a notable interest in the
mathematical literature (see discussion below), 
seem not to be so far studied in the case of 
logarithmic potentials.

The sixth order system \eqref{CH1}-\eqref{CH2} arises as a model
of dynamics of ternary oil-water-surfactant mixtures in which three
phases occupying a region
$\Omega$ in $ \RR^3$, microemulsion, almost pure oil and almost pure water, 
can coexist in equilibrium.
The phenomenological Landau-Ginzburg theory for such mixtures has been proposed
in a series of papers by Gompper et~al.~(see, e.g., 
\cite{GK93a,GK93b,GZ92} and other references in \cite{PZ11}).
This theory is based on the free energy functional
\eqref{defiE} with constant $\delta>0$ (in general, however,
this coefficient can depend on $u$, see \cite{SS93}), and with
$F(u)$, $a(u)$ approximated, respectively, by a sixth and
a second order polynomial:
\begin{equation}\label{apprFa}
  F(u)= (u+1)^2 (u^2+h_0) (u-1)^2, \qquad
   a(u) = g_0 + g_2 u^2,
\end{equation}
where the constant parameters $h_0,g_0,g_2$ are adjusted experimentally, 
$g_2>0$ and $h_0$, $g_0$ are of arbitrary sign. 
In this model, $u$ is the scalar, conserved order parameter 
representing the local difference
between oil and water concentrations; $u=-1$, $u=1$, and $u=0$ correspond to oil-rich, 
water-rich and microemulsion phases, respectively, and the parameter $h_0$ measures 
the deviation from oil-water-microemulsion coexistence. 

The associated evolution system \eqref{CH1}-\eqref{CH2} has the standard
Cahn-Hilliard structure. Equation \eqref{CH1} expresses the conservation law
\begin{equation}\label{conslaw}  
      u_t + \nabla \cdot j = 0
\end{equation}
with the mass flux $j$ given by
\begin{equation}\label{mflux}
j = - M\nabla w.
\end{equation}
Here $M > 0$ is the constant mobility (we set $M=1$ for simplicity),
and $w$ is the chemical potential difference between the oil and water phases. 
The chemical potential is defined by the constitutive equation
\begin{equation}\label{chpot}  
 w = \frac{\delta \calE\dd(u)}{\delta u} + \epsi u_t,
\end{equation}
where $ \frac{\delta \calE\dd(u)}{\delta u}$ is the first variation of the functional
$ \calE\dd(u) $, and the constant $ \epsi \geq 0$
represents possible viscous effects.
For energy \eqref{defiE} equation \eqref{chpot} yields \eqref{CH2}.
We note also that the boundary conditions $ \dn u =  \delta \dn \Delta u = 0 $
are standardly used in the frame of sixth order Cahn-Hilliard
models due to their mathematical simplicity. Moreover,
they are related to the variational structure of the 
problem in terms of the functional \eqref{defiE}.
However, other types of boundary conditions for $u$ might 
be considered as well, paying the price of technical 
complications in the proofs. Concerning, instead, 
the condition $ \dn w = 0$, in view of \eqref{mflux}, 
it simply represents the mass isolation at the 
boundary of $\Omega$.

The  system \eqref{CH1}-\eqref{neum-intro}
with functions $ F(u), a(u)$ in the polynomial form \eqref{apprFa},
and with no viscous term
($\epsi=0$) has been recently 
studied in~\cite{PZ11}. It has been proved there that for a sufficiently 
smooth initial datum $u_0$ the system admits a unique global solution in the strong sense.

The sixth order Cahn-Hilliard type equation with the same structure as 
\eqref{CH1}-\eqref{CH2}, $\delta > 0$, polynomial $F(u)$, and negative constant $a$,
arises also as the so-called phase field crystal (PFC)
atomistic model of crystal growth, developed by Elder et al.,
see e.g., \cite{EG04, BGE06, BEG08}, and \cite{GDL09} for the overview and up-to 
date references. It is also worth mentioning
a class of sixth order convective Cahn-Hilliard type equations 
with different (nonconservative) structure than \eqref{CH1}-\eqref{CH2}.
This type of equations arise in particular as a model of
the faceting a of a growing crystalline surface, 
derived by Savina et al.~\cite{S03}
(for a review of other convective 4th and 6th order Cahn-Hilliard models see \cite{KEMW08}). 
In this model, contrary to~\eqref{CH1}-\eqref{CH2}, the 
order parameter $u$ is not a conserved quantity due to the 
presence of a force-like term related to the deposition rate.
Such class of models has been recently studied mathematically
in one- and two- dimensional cases by Korzec et al.~\cite{KEMW08, KNR11,KR11}.

Finally, let us note that in the case $\delta=0$, $a(u)=\const>0$, 
the functional \eqref{defiE}
represents the classical Cahn-Hilliard free energy \cite{Ca,CH}. 
The original Cahn-Hilliard free energy derivation
has been extended by Lass et al. \cite{LJS06}
to account for composition dependence of the
gradient energy coefficient $a(u)$. For a face-centered cubic
crystal the following expressions for $a(u)$ have been derived,
depending on the level of approximation of the nearest-neighbor
interactions:
\begin{equation}\label{appra}
   a(u) = a_0 + a_1 u + a_2 u^2,
\end{equation}
where $a_0>0$, $a_1,a_2\in\RR$ in the case of four-body
interactions, $a_2=0$ in the case of three-body
interactions, and $a_1=a_2=0$ in the case of pairwise
interactions. 
Numerical experiments in \cite{LJS06} indicate that these
three different approximations (all reflecting the face-centered 
cubic crystal symmetry) have a substantial effect on the shape 
of the equilibrium composition profile and the
interfacial energy. 

A specific free energy with composition dependent gradient
energy coefficient $a(u)$ also arises in modelling of phase
separation in polymers \cite{dG80}. This energy, known
as the Flory-Huggins-de Gennes one, has the form \eqref{defiE} 
with $\delta=0$, $F(u)$ being the logarithmic potential
\eqref{logpot}, and the singular coefficient
\begin{equation}\label{adG}
   a(u) = \frac1{(1-u)(1+u)}.
\end{equation}
We mention also that various formulations of phase-field models
with gradient energy coefficient dependent on the order parameter
(and possibly on other fields) appear, e.g., in~\cite{Aif86,BS96}.

Our objective in this paper is threefold. 
First, we would like to extend the result 
of \cite{PZ11} both to the viscous problem 
($\epsi>0$) and to the case 
when the configuration potential is {\sl singular}\ 
(e.g., of the form \eqref{logpot}).
While the first extension is almost straighforward, considering
constraint (singular) terms in fourth order equations 
(\eqref{CH2}, in the specific case) gives rise 
to regularity problems since it is not possible, up to
our knowledge, to estimate all the terms of equation~\eqref{CH2}
in $L^p$-spaces. For this reason, the nonlinear term $f(u)$ 
has to be intended in a weaker form, namely, as a selection
of a nonlinear, and possibly multivalued, mapping acting
from $V=H^1(\Omega)$ to $V'$. This involves some monotone 
operator technique that is developped in a specific section
of the paper.

As a second step, we investigate the behavior of the solutions
to the sixth order system as the parameter $\delta$ is let to tend
to $0$. In particular, we would like to show that, at least up to
subsequences, we can obtain in the limit suitably defined
solutions to the fourth order system obtained setting $\delta = 0$ 
in \eqref{CH2}. Unfortunately, we are able to prove this fact 
only under additional conditions. The reason is that
the natural estimate required to control second space derivatives
of $u$, i.e., testing \eqref{CH2} by $-\Delta u$, is compatible with 
the nonlinear term in $\nabla u$ only under additional assumptions 
on $a$ (e.g., if $a$ is concave). This nontrivial fact depends
on an integration by parts formula devised by Dal Passo, 
Garcke and Gr\"un in \cite{DpGG} in the frame 
of the thin-film equation and whose use is necessary to control
the nonlinear gradient term. It is however likely that the use of more
refined integration by parts techniques may permit to control
the nonlinear gradient term under more general conditions on $a$.

Since we are able to take the limit $\delta\searrow0$ only in special
cases, in the subsequent part of the paper we address
the fourth order problem by using a direct approach. 
In this way, we can obtain existence of a weak solution 
under general conditions on $a$ (we notice that, however, 
uniqueness is no longer guaranteed for $\delta=0$).
The proof of existence is based on an ``ad hoc''
regularization of the equations by means of a system of
phase-field type. This kind of approach has been proved
to be effective also in the frame of other types
of Cahn-Hilliard equations (see, e.g., \cite{BaPa05}).
Local existence for the regularized system
is then shown by means of the Schauder theorem, and, finally,
the regularization is removed by means of suitable a priori estimates
and compactness methods. This procedure involves 
some technicalities since parabolic spaces of H\"older type
have to be used for the fixed point argument.
Indeed, the use of Sobolev techniques seems not suitable due to 
the nonlinearity in the highest order term, which prevents
from having compactness of the fixed point map with respect to
Sobolev norms. A further difficulty is related with the
necessity of estimating the second order space derivatives
of $u$ in presence of the nonlinear term in the gradient. This is
obtained by introducing a proper transformed variable, 
and rewriting \eqref{CH2} in terms of it. Proceeding in 
this way, we can get rid of that nonlinearity, but
at the same time, we can still
exploit the good monotonicity properties of $f$. 
We note here that a different method based on entropy estimates
could also be used to estimate $\Delta u$ without making the change
of variable, which seems however a simpler technique.

Finally, in the last section of the paper, we 
discuss further property
of weak solutions. More precisely, we address the problems 
of uniqueness (only for the 4th order system,
since in the case $\delta>0$ it is always guaranteed)
and of parabolic time-regularization of solutions
(both for the 6th and for the 4th order system).
We are able to prove such properties only when the
energy functional $\calE\dd$ is $\lambda$-convex
(so that its gradient is monotone up to a linear
perturbation). In terms of the coefficient $a$,
this corresponds to asking that $a$ is a {\sl convex}\/ 
function and, moreover, $1/a$ is {\sl concave}\/
(cf.~\cite{DNS} for generalizations and further comments
regarding this condition). If these conditions
fail, then the gradient of the energy functional exhibits
a nonmonotone structure in terms of the space derivatives
of the highest order. For this reason, proving an
estimate of contractive type (which would be required for
having uniqueness) appears to be difficult in that case.

As a final result, we will show that, both in the 
6th and in the {\sl viscous}\/ 4th order case, 
all weak solutions satisfy the energy {\sl equality}\/
\eqref{energyineq}, at least in an integrated form 
(and not just an energy inequality).
This property is the starting point for proving
existence of the global attractor for 
the dynamical process associated
to system \eqref{CH1}-\eqref{CH2}, an issue that 
we intend to investigate in a forthcoming paper.
Actually, it is not difficult to show that
the set of initial data having finite energy
constitutes a complete metric space
(see, e.g., \cite[Lemma~3.8]{RS}) which can
be used as a {\sl phase space}\/ for the system. 
Then, by applying the so-called ``energy method''
(cf., e.g., \cite{MRW,Ba1}), one can see that
the energy equality implies precompactness
of trajectories for $t\nearrow\infty$ with
respect to the metric of the phase space. In
turn, this gives existence of the global attractor
with respect to the same metric.
On the other hand, the question whether the energy equality 
holds in the nonviscous 4th order case seems
to be more delicate, and, actually, we could not 
give a positive answer to it.

It is also worth to notice an important issue concerned
with the sharp interface limit of the Cahn-Hilliard equation 
with a nonlinear gradient energy coefficient $a(u)$.
To our knowledge this issue has not been so far addressed in 
the literature. Let us mention that using the method of matched
asymptotic expansions the sharp interface limits of the 
Cahn-Hilliard equation with constant coefficient $a$ have been
investigated by Pego \cite{Peg89} and rigorously by Alikakos et al.
\cite {ABC94}.
Such method has been also successfully applied to a number of 
phase field models of phase transition problems, see e.g., \cite{CF88},
\cite{C90}. 
In view of various physical applications described above, it would 
be of interest to apply the matched asymptotic expansions in the case 
of a nonlinear coefficient $a(u)$ to investigate what kind of corrections it 
may introduce to conditions on the sharp interface.

The plan of the paper is as follows. 
In the next Section~\ref{sec:main},
we will report our notation and hypotheses, together with
some general tools that will be used in the proofs. 
Section~\ref{sec:6th} will contain
the analysis of the sixth order model. The limit 
$\delta\searrow 0$ will then be analyzed in 
Section~\ref{sec:6thto4th}.
Section~\ref{sec:4th} will be devoted to the analysis 
of the fourth order model. Finally, in Section~\ref{sec:uniq}
uniqueness and regularization properties of the solutions
will be discussed, as well as the validity of the
energy equality.

\medskip

\noindent%
{\bf Acknowledgment.}~~The authors are grateful to Prof.~Giuseppe 
Savar\'e for fruitful discussions about the strategy of some proofs.


\section{Notations and technical tools}
\label{sec:main}

Let $\Omega$ be a smooth bounded domain of $\RR^3$
of boundary $\Gamma$, $T>0$ a given final time, and
let $Q:=(0,T)\times\Omega$. We let $H:=L^2(\Omega)$,
endowed with the standard scalar product $(\cdot,\cdot)$
and norm $\| \cdot \|$. 
For $s>0$ and $p\in[1,\infty]$, we 
use the notation $W^{s,p}(\Omega)$ to indicate 
Sobolev spaces of positive (possibly fractional) order.
We also set $H^s(\Omega):=W^{s,2}(\Omega)$ and
let $V:=H^1(\Omega)$.
We note by $\duav{\cdot,\cdot}$ the duality between
$V'$ and $V$ and by $\|\cdot\|_X$ the norm in 
the generic Banach space $X$. 
We identify $H$ with
$H'$ in such a way that $H$ can be seen as a subspace of $V'$
or, in other words, $(V,H,V')$ form a Hilbert triplet.

We make the following assumptions on the nonlinear
terms in \eqref{CH1}-\eqref{CH2}:
\begin{align}\label{hpa1}
  & a \in C^2_b(\RR;\RR), \quad \esiste \agiu,\asu>0:~~
    \agiu \le a(r)\le \asu~~\perogni r\in \RR;\\
 \label{hpa2}
  & \esiste a_-,a_+\in [\agiu,\asu]:~~
   a(r)\equiv a_-~~\perogni r\le-2, \quad
   a(r)\equiv a_+~~\perogni r\ge 2;\\
 \label{hpf1}
  & f\in C^1((-1,1);\RR), \quad f(0)=0, \quad
   \esiste\lambda\ge 0:~~f'(r)\ge -\lambda~~\perogni r\in (-1,1);\\
 \label{hpf2}
  & \lim_{|r|\to 1}f(r)r = \lim_{|r|\to 1}\frac{f'(r)}{|f(r)|}
   = + \infty.
\end{align}
The latter condition in \eqref{hpf2} is just a technical hypotheses
which is actually verified in all significant cases.
We also notice that, due to the choice of a singular potential
(mathematically represented here by assumptions \eqref{hpf1}-\eqref{hpf2}),
any weak solution $u$ will take its values only in the physical 
interval $[-1,1]$). For this reason, the behavior of 
$a$ is also significant only in that interval and we have extended it
outside $[-1,1]$ just for the purpose of properly constructing
the approximating problem (see Subsection~\ref{subsec:appr} below). 

Note that our assumptions on $a$ are not in conflict with 
\eqref{apprFa} or \eqref{appra} since these conditions
(or, more generally, any condition on the values of $a(u)$ 
for large $u$) make sense in the different situation of a
function $f$ with polynomial growth (which does not constrain
$u$ in the interval $(-1,1)$).
It should be pointed out, however, that the assumptions \eqref{hpa1}-\eqref{hpf2}
do not admit the singular Flory-Huggins-de Gennes free energy model
with $a(u)$ given by \eqref{adG}. We expect that the analysis of such 
a singular model could require different techniques.

In \eqref{hpa1}, $C^2_b$ denotes
the space of functions that are continuous and globally
bounded together with their derivatives up to the second order.
Concerning $f$, \eqref{hpf1} states that it can be written in the form
\begin{equation}\label{f0}
  f(r)=f_0(r)-\lambda r,
\end{equation}
i.e., as the difference between a (dominating) monotone part $f_0$ and
a linear perturbation. By \eqref{hpf1}-\eqref{hpf2}, we can also set,
for $r\in(-1,1)$,
\begin{equation}\label{F0}
  F_0(r):=\int_0^r f_0(s)\,\dis \qquext{and~}\,
   F(r):=F_0(r)-\frac\lambda2 r^2,
\end{equation}
so that $F'=f$.
Notice that $F_0$ may be bounded in $(-1,1)$ (e.g., this occurs
in the case of the logarithmic potential \eqref{logpot}). If this 
is the case, we extend it by continuity to $[-1,1]$. Then, 
$F_0$ is set to be $+\infty$ either outside $(-1,1)$ (if it is 
unbounded in $(-1,1)$) or outside $[-1,1]$ (if it is bounded
in $(-1,1)$). This standard procedure permits to penalize the 
non-physical values of the variable $u$ and to intend $f_0$ as the
subdifferential of the (extended)
convex function $F_0:\RR\to[0,+\infty]$.

That said, we define a number of operators. First, we set
\begin{equation}\label{defiA}
   A:V\to V', \qquad
    \duav{A v, z}:= \io \nabla v \cdot \nabla z,
     \quext{for }\, v,z \in V.
\end{equation}
Then, we define 
\begin{equation}\label{defiW}
  W:=\big\{z\in H^2(\Omega):~\dn z=0~\text{on }\Gamma\big\}
\end{equation}
and recall that (a suitable restriction of) $A$ can be seen
as an unbounded linear operator on $H$
having domain $W$. The space
$W$ is endowed with the natural $H^2$-norm. We then
introduce
\begin{equation}\label{deficalA}
   \calA: W \to H, \qquad
    \calA(z) := - a(z)\Delta z - \frac{a'(z)}2 |\nabla z|^2.
\end{equation}
It is a standard issue to check that, indeed, $\calA$ takes its
values in~$H$.


\subsection{Weak subdifferential operators}
\label{sec:weak}

To state the weak formulation of the 6th order system, 
we need to introduce a proper 
relaxed form of the maximal monotone
operator associated to the function $f_0$ and 
acting in the duality between $V'$ and $V$ (rather than
in the scalar product of $H$). Actually, it is well known
(see, e.g., \cite[Ex.~2.1.3, p.~21]{Br})
that $f_0$ can be interpreted as a maximal monotone operator
on~$H$ by setting, for $v,\xi\in H$,
\begin{equation} \label{betaL2}
  \xi = f_0(v)\quext{in $H$}~~~\Longleftrightarrow~~~
   \xi(x) = f_0(v(x))\quext{a.e.~in $\Omega$}.
\end{equation}
If no danger of confusion occurs,
the new operator on $H$ will be still noted by the letter $f_0$.
Correspondingly, $f_0$ is the $H$-subdifferential of the 
convex functional
\begin{equation} \label{betaL2-2}
  \calF_0:H\mapsto[0,+\infty], \qquad
   \calF_0(v):= \io F_0(v(x)),   
\end{equation}
where the integral might possibly be $+\infty$ (this happens,
e.g., when $|v|>1$ on a set of strictly positive Lebesgue measure).

The weak form of $f_0$ can be introduced by setting
\begin{equation} \label{betaV}
  \xi\in \fzw(v) \Longleftrightarrow
   \duav{\xi,z-v}\le \calF_0(z)-\calF_0(v)
    \quext{for any $z\in V$}. 
\end{equation}
Actually, this is nothing else than the 
definition of the subdifferential
of (the restriction to $V$ of) $\calF_0$ 
with respect to the duality pairing
between $V'$ and $V$. In general, $\fzw$ can be a 
{\sl multivalued}\/ operator; namely,
$\fzw$ is a {\sl subset}\ of $V'$ that
may contain more than one element.
It is not difficult to prove 
(see, e.g., \cite[Prop.~2.5]{BCGG}) that, if $v\in V$
and $f_0(v)\in H$, then
\begin{equation} \label{betavsbetaw}
  \{f_0(v)\}\subset\fzw(v).
\end{equation}
Moreover, 
\begin{equation} \label{betavsbetaw2}
  \text{if }\,v\in V~\,\text{and }\,
  \xi \in \fzw(v) \cap H, 
   \quext{then }\,\xi = f_0(v)
   ~\,\text{a.e.~in }\,\Omega.
\end{equation}
In general, the inclusion in \eqref{betavsbetaw}
is strict and, for instance, it can happen that 
$f_0(v)\not\in H$ (i.e., $v$ does not belong
to the $H$-domain of $f_0$), 
while $\fzw(v)$ is nonempty. Nevertheless,
we still have some ``automatic'' gain of 
regularity for any element of $\fzw(v)$:
\bepr\label{misura}
 Let $v\in V$, $\xi\in\fzw(v)$. Then,
 $\xi$ can be seen as an element of the 
 space ${\cal M}({\overline \Omega})=C^0(\barO)'$
 of the bounded real-valued Borel measures on $\overline \Omega$. 
 More precisely, there exists $T\in {\cal M}({\overline \Omega})$,
 such that 
 \begin{equation} \label{identif}
   \duav{\xi,z}=\ibaro z\,\diT
    \qquext{for any~\,$z\in V\cap C^0(\overline\Omega)$}.
 \end{equation}
\empr
\begin{proof}
 Let $z\in C^0(\overline \Omega)\cap V$ with $z\not = 0$.
 Then, using definition \eqref{betaV},
 it is easy to see that 
 \begin{align} \label{prova-meas}
   \duav{\xi,z} & = 2\| z \|_{L^\infty(\Omega)} \duavg{\xi,\frac{z}{2\| z \|_{L^\infty(\Omega)}}}
    \le 2\| z \|_{L^\infty(\Omega)} \bigg(
       \duav{\xi,v}+\calF_0\Big( \frac{z}{2\| z \|_{L^\infty(\Omega)}} \Big)-\calF_0(v) \bigg)\\
     &  \le 2\| z \|_{L^\infty(\Omega)} \Big( |\duav{\xi,v}|+|\Omega|\big(F_0(-1/2)+F_0(1/2)\big) \Big).
 \end{align}
 This actually shows that the linear 
 functional $z\mapsto\duav{\xi,z}$
 defined on $C^0(\overline \Omega)\cap V$ (that is a dense subspace
 of $C^0(\overline \Omega)$, recall that $\Omega$ is smooth) 
 is continuous with respect to the sup-norm. Thus, 
 by the Riesz representation theorem, it can 
 be represented over $C^0(\barO)$
 by a measure $T\in {\cal M}({\overline \Omega})$.
\end{proof}

\noindent%
Actually, we can give a general definition, saying that a functional 
$\xi\in V'$ belongs to the space 
$V'\cap {\cal M}(\overline \Omega)$ provided that 
$\xi$ is continuous 
with respect to the sup-norm on $\overline \Omega$.
In this case, we can use \eqref{identif} and say
that the measure $T$ represents $\xi$ on 
${\cal M}(\overline \Omega)$. 
We now recall a result \cite[Thm.~3]{brezisart}
that will be exploited in the sequel.
\bete\label{teobrezis}
 Let $v\in V$, $\xi\in \fzw(v)$. Then, 
 denoting by $\xi_a+\xi_s=\xi$ the Lebesgue decomposition 
 of $\xi$, with $\xi_a$ ($\xi_s$) standing for
 the absolute continuous (singular, respectively)
 part of $\xi$, we have
 \begin{align}\label{bre1}
   & \xi_a v\in L^1(\Omega),\\
  \label{bre2}
   & \xi_a(x) = f_0(v(x)) \qquext{for a.e.~$x\in\Omega$,}\\
  \label{bre3}
   & \duav{\xi,v} - \io \xi_a v\,\dix
    = \sup \bigg\{\ibaro z\,\dixi_s,~z\in C^0(\overline\Omega),~
    z(\overline\Omega)\subset[-1,1] \bigg\}.
 \end{align}
\ente
\noindent%
Actually, in \cite{brezisart} a slightly different result is proved, 
where $V$ is replaced by $H^1_0(\Omega)$ and, correspondingly,
${\cal M}(\overline \Omega)$ is replaced by ${\cal M}(\Omega)$
(i.e., the dual of $C_c^0(\Omega)$). Nevertheless, thanks to the
smoothness of $\Omega$, one can easily realize that the approximation
procedure used in the proof of the theorem can be extended to 
cover the present situation. The only difference is given by
the fact that the singular part $\xi_s$ may be supported
also on the boundary.

\smallskip

Let us now recall that, given a pair $X,Y$ of Banach spaces, 
a sequence of (multivalued)
operators ${\cal T}_n:X\to 2^Y$ is said to 
G-converge (strongly) to ${\cal T}$ iff
\begin{equation}\label{defGconv}
  \perogni (x,y)\in {\cal T}, \quad \esiste 
   (x_n,y_n)\in {\cal T}_n \quext{such that \,
   $(x_n,y_n)\to(x,y)$~~strongly in }\, X\times Y.
\end{equation}
We would like to apply this condition to an approximation
of the monotone function $f_0$ that we now construct. 
Namely, for $\sigma\in(0,1)$ (intended
to go to $0$ in the limit), we would like to
have a family $\{f\ssi\}$ of monotone functions such that
\begin{align}\label{defifsigma}
  & f\ssi\in C^{1}(\RR), \qquad
   f\ssi'\in L^{\infty}(\RR), \qquad
   f\ssi(0)=0, \\
 \label{convcomp}  
  & f\ssi\to f_0 \quext{uniformly on compact
   subsets of }\,(-1,1).
\end{align}
Moreover, noting 
\begin{equation}\label{defiFsigma}
  F\ssi(r):=\int_0^r 
   f\ssi(s)\,\dis,
   \quext{for }\,r\in\RR,
\end{equation}
we ask that
\begin{equation}\label{propFsigma}
  F\ssi(r) \ge \lambda r^2 - c,
\end{equation}
for some $c\ge 0$ independent of $\sigma$ and for
all $r\in\RR$, $\sigma\in (0,1)$,
where $\lambda$ is as in
\eqref{hpf1} (note that the analogue
of the above property holds for $F$ thanks to 
the first of \eqref{hpf2}). 
Moreover, we ask the monotonicity
condition
\begin{equation}\label{defifsigma2}
  F_{\sigma_1}(r)\le F_{\sigma_2}(r)
   \qquext{if }\,\sigma_2\le \sigma_1
   \quext{and for all }\, r\in \RR.
\end{equation}
Finally, on account of the last assumption
\eqref{hpf2}, we require that
\begin{equation}\label{goodmono}
  \perogni m>0,~~\esiste C_m\ge 0:~~~
   f\ssi'(r) - m |f\ssi(r)| \ge - C_m,
   \quad \perogni r\in[-2,2]
\end{equation}
with $C_m$ being independent of $\sigma$.
Notice that it is sufficient to ask 
the above property for $r\in [-2,2]$. 
The details of the construction of a family $\{f\ssi\}$ 
fulfilling \eqref{defifsigma}-\eqref{goodmono} are 
standard and hence we leave them to the reader.
For instance, one can first take Yosida regularizations
(see, e.g., \cite[Chap.~2]{Br}) and then mollify in order
to get additional smoothness.

Thanks to the monotonicity property \eqref{defifsigma2},
we can apply \cite[Thm.~3.20]{At}, which gives that
\begin{align}\label{Gforte}
  & f\ssi\quext{G-converges to }\,f_0
   \quext{in \,$H\times H$},\\
 \label{Gdebole}
  & f\ssi\quext{G-converges to }\,\fzw
   \quext{in \,$V\times V'$}.
\end{align}
A notable consequence of G-convergence is
the following property, whose proof can be obtained by
slightly modifying \cite[Prop.~1.1, p.~42]{barbu}:
\bele\label{limimono}
 Let $X$ be an Hilbert space, ${\cal B}\ssi$, ${\cal B}$ be 
 maximal monotone operators in $X\times X'$ such that 
 \begin{equation}\label{Gastr}
   {\cal B}\ssi\quext{G-converges to }\, {\cal B}
    \quext{in }\,X\times X',
 \end{equation}
 as $\sigma\searrow0$. Let also, for any $\sigma>0$, $v\ssi\in X$, 
 $\xi\ssi\in X'$ such that $\xi\ssi\in{\cal B}\ssi(v\ssi)$.
 Finally, let us assume that, for some $v\in X$, 
 $\xi\in X'$, there holds
 \begin{align}\label{Gastrnew}
   & v\ssi\to v\quext{weakly in }\, X, \qquad
    \xi\ssi\to \xi\quext{weakly in }\,X',\\[1mm]
  \label{Gastrnew-2}
   & \limsup_{\sigma\searrow0} \duavg{\xi\ssi,v\ssi}_X
    \le \duavg{\xi,v}_X.
 \end{align}
 Then, $\xi\in {\cal B}(v)$.
\enle
\noindent%
Next, we present an integration by parts formula:
\bele\label{BSesteso}
 Let $u\in W\cap H^3(\Omega)$, $\xi\in V'$ such that
 $\xi\in \fzw(u)$. Then, we have that
 \begin{equation}\label{majozero}
   \duav{\xi,Au}\geq 0.
 \end{equation}
\enle
\begin{proof}
Let us first note that the duality above surely makes
sense in the assigned regularity setting. Actually,
we have that $Au\in V$. We then consider the elliptic problem
\begin{equation}\label{elpromon}
  u\ssi\in V, \qquad
   u\ssi + A^2 u\ssi + f\ssi (u\ssi)
    = u + A^2 u + \xi \text{~~~~in \,$V'$.}
\end{equation}
Since $f\ssi$ is monotone and Lipschitz continuous
and the above \rhs\ lies in $V'$,
it is not difficult to show that the above problem
admits a unique solution $u\ssi \in W \cap H^3(\Omega)$.

Moreover, the standard a priori estimates for $u\ssi$
lead to the following convergence relations, which hold, 
for some $v\in V$ and $\zeta\in V'$,
up to the extraction of (non-relabelled)
subsequences (in fact uniqueness
guarantees them for the whole $\sigma\searrow0$):
\begin{align}\label{stlemma11}
  & u\ssi\longrightarrow v \quext{weakly in }\, H^3(\Omega)~~
   \text{and strongly in }\,W,\\
 \label{stlemma11.2}
  &  A^2 u\ssi\longrightarrow A^2 v \quext{weakly in }\, V',\\
 \label{stlemma11.3}
  &  f\ssi(u\ssi)\longrightarrow \zeta \quext{weakly in }\, V'.
\end{align}
As a byproduct, the limit functions satisfy
$ v+A^2v+\zeta=u+A^2u+\xi$ in~$V'$. Moreover, we
deduce from (\ref{elpromon})
\begin{equation}\label{contolemma11}
  \big(f\ssi(u\ssi),u\ssi\big)
   = \duavg{ u + A^2 u + \xi - u\ssi - A^2u\ssi, u\ssi},
\end{equation} 
whence 
\begin{equation}\label{old11}
  \lim_{\sigma \rightarrow 0} \,\big(f\ssi(u\ssi),u\ssi\big)
   = \duavg{u + A^2 u + \xi - v - A^2 v, v}
   = \duav{\zeta,v}.
\end{equation}
Then, on account of 
\eqref{stlemma11}, \eqref{stlemma11.3}, \eqref{Gdebole}
and Lemma~\ref{limimono} (cf., in particular,
relation \eqref{Gdebole}) applied to the sequence
$\{f_\sigma(v_\sigma)\}$, we readily obtain that 
$\zeta\in \fzw(v)$. By uniqueness, $v=u$ and $\zeta=\xi$. 

Let us finally verify the required property.
Actually, for $\sigma>0$, thanks to monotonicity
of $f\ssi$ we have
\begin{equation}\label{contolemma12}
  0 \leq \big(f\ssi(u\ssi), A u\ssi \big)
   = \duavg{ u + A^2 u+\xi-u\ssi-A^2u\ssi, A u\ssi}.
\end{equation}
Taking the supremum limit, we then obtain
\begin{equation}\label{contolemma12b}
  0 \leq \limsup_{\sigma\searrow 0} \duavg{ u + A^2 u+\xi-u\ssi-A^2u\ssi, A u\ssi}
   = \duavg{ u + A^2 u+ \xi - u, Au}  - \liminf_{\sigma\searrow 0} \duavg{A^2u\ssi, A u\ssi}.
\end{equation}
Then, using \eqref{stlemma11} and semicontinuity of norms 
with respect to weak convergence, 
\begin{equation}\label{contolemma12c}
  - \liminf_{\sigma\searrow 0} \duavg{A^2u\ssi, A u\ssi}
   = - \liminf_{\sigma\searrow 0} \| \nabla A u\ssi \|^2
   \le - \| \nabla A u \|^2
   = - \duavg{A^2 u , A u},
\end{equation}
whence we finally obtain
\begin{equation}\label{old12} 
  0 \le \duav{u+A^2u+\xi-u-A^2u,Au}
   = \duav{\xi,Au},
\end{equation}
as desired.
\end{proof}
\noindent%
Next, we recall a further integration by parts formula
that extends the classical result 
\cite[Lemma~3.3, p.~73]{Br}
(see, e.g., \cite[Lemma~4.1]{RS} for a proof):
\bele\label{BResteso}
 Let $T>0$ and let $\calJ:H\to [0,+\infty]$ a convex,
 lower semicontinuous and proper functional. Let 
 $u\in \HUVp \cap \LDV$, $\eta\in \LDV$ and let 
 $\eta(t)\in \de\calJ(u(t))$ for a.e.~$t\in(0,T)$,
 where $\de\calJ$ 
 is the $H$-subdifferential of $\calJ$. Moreover, let us
 suppose the coercivity property
 \begin{equation}\label{coerccalJ}
   \esiste k_1>0,~k_2\ge 0
   \quext{such that }\,\calJ(v) \ge k_1 \| v \|^2 - k_2
   \quad\perogni v\in H.
 \end{equation}
 Then, the function $t\mapsto \calJ(u(t))$ is
 absolutely continuous in $[0,T]$ and 
 \begin{equation}\label{ipepardiff}
   \ddt \calJ(u(t)) = \duav{u_t(t),\eta(t)}
    \quext{for a.e.~}\, t\in (0,T).
 \end{equation}
 In particular, integrating in time, we have
 \begin{equation}\label{ipepars}
   \int_s^t \duav{u_t(r),\eta(r)}\,\dir
    = \calJ(u(t)) - \calJ(u(s))
    \quad\perogni s,t\in [0,T].
 \end{equation}
\enle
\noindent%
We conclude this section by stating an integration by parts 
formula for the operator $\calA$.
\bele\label{lemma:ipp}
 Let $a$ satisfy \eqref{hpa1} and let either
 \begin{equation} \label{x11}
   v \in \HUH \cap \LDW \cap L^\infty(Q),
 \end{equation}  
 or
 \begin{equation} \label{x12}
   v \in \HUVp \cap \LIW \cap L^2(0,T;H^3(\Omega)).
 \end{equation}  
 Then, the function
 \begin{equation} \label{x13}
   t\mapsto \io \frac{a(v(t))}2 | \nabla v(t) |^2
 \end{equation}  
 is absolutely continuous over $[0,T]$. Moreover, 
 for all $s,t\in [0,T]$ we have that
 \begin{equation} \label{x14}
   \int_s^t \big(\calA(v(r)),v_t(r)\big)\,\dir
    = \io \frac{a(v(t))}2 |\nabla v(t)|^2
    - \io \frac{a(v(s))}2 |\nabla v(s)|^2,
 \end{equation}  
 where, in the case \eqref{x12}, the scalar product 
 in the integral on 
 the \lhs\ has to be replaced with the duality 
 $\duav{v_t(r),\calA(v(r))}$.
\enle
\begin{proof}
We first notice that \eqref{x13}-\eqref{x14}
surely hold if $v$ is smoother. Then, we can proceed
by first regularizing $v$ and then passing to the
limit. Namely, we define $v\ssi$, a.e.~in~$(0,T)$,
as the solution of the singular perturbation problem
\begin{equation} \label{co93}
  v\ssi + \sigma A v\ssi = v,
   \quext{for }\, \sigma\in (0,1).
\end{equation}  
Then, in the case \eqref{x11}, we have 
\begin{equation} \label{co93-b}
  v\ssi \in H^1(0,T;W) \cap L^2(0,T;H^4(\Omega)),
\end{equation}  
whereas, if \eqref{x12} holds, we get
\begin{equation} \label{co93-c}
  v\ssi \in H^1(0,T;V) \cap L^\infty(0,T;H^4(\Omega)).
\end{equation}  
Moreover, proceeding as in \cite[Appendix]{CGG}
(cf., in particular, Proposition~6.1 therein)
and applying the Lebesgue dominated convergence
theorem in order to control the dependence on time variable,
we can easily prove that
\begin{equation} \label{y11}
  v\ssi \to v \quext{strongly in }\,\HUH \cap \LDW   
   ~~\text{and weakly star in }\, L^\infty(Q)
\end{equation}  
(the latter condition following from the maximum principle),
if \eqref{x11} holds, or 
\begin{equation} \label{y12}
  v\ssi \to v \quext{strongly in }\,\HUVp \cap L^2(0,T;H^3(\Omega)) 
   ~~\text{and weakly star in }\, \LIW,
\end{equation}  
if \eqref{x12} is satisfied instead.

Now, the functions $v\ssi$, being smooth,
surely satisfy the analogue of \eqref{x14}:
\begin{equation} \label{x14ssi}
  \int_s^t \big(\calA(v\ssi(r)),v_{\sigma,t}(r)\big)\,\dir
   = \io \frac{a(v\ssi(t))}2 |\nabla v\ssi(t)|^2
   - \io \frac{a(v\ssi(s))}2 |\nabla v\ssi(s)|^2,
\end{equation}  
for all $s,t\in[0,T]$.
Let us prove that we can take the limit $\sigma\searrow0$,
considering first the case \eqref{x11}. Then, using \eqref{y11}
and standard compactness results, it is not difficult 
to check that
\begin{equation} \label{co93-b2}
  \calA(v\ssi) \to \calA(v), \quext{(at least) weakly in }\,L^2(0,T;H).
\end{equation}  
In particular, to control the square gradient term in $\calA$, 
we use the Gagliardo-Nirenberg inequality (cf.~\cite{Ni})
\begin{equation}\label{ineq:gn}
  \| \nabla z \|_{L^4(\Omega)} \le c\OO \| z \|_{W}^{1/2}
   \| z \|_{L^\infty(\Omega)}^{1/2}
   + \| z \| \qquad \perogni z \in W,
\end{equation}
so that, thanks also to \eqref{hpa1},
\begin{equation}\label{conseq:gn}
  \big\| a'(v\ssi) |\nabla v\ssi|^2 \big\|_{\LDH}
   \le \| a'(v\ssi) \|_{L^\infty(Q)} \| \nabla v\ssi \|_{L^4(Q)}^2
   \le c \| v\ssi \|_{L^\infty(Q)} \big( 1 +
     \| A v\ssi \|_{\LDH} \big),
\end{equation}
and \eqref{co93-b2} follows.
Moreover, by \eqref{y11} and the continuous embedding
$H^1(0,T;H) \cap L^2(0,T;W) \subset C^0([0,T];V)$,
we also have that
\begin{equation} \label{co93-d}
  v\ssi \to v \quext{strongly in }\,C^0([0,T];V).
\end{equation}  
Combining \eqref{y11}, \eqref{co93-b2} and \eqref{co93-d}, we
can take the limit $\sigma\searrow 0$ in \eqref{x14ssi} and get
back \eqref{x14}. Then, the absolute continuity property
of the functional in \eqref{x13} follows from the summability
of the integrand on the \lhs\ of \eqref{x14}.

Finally, let us come to the case \eqref{x12}. Then,
\eqref{y12} and the Aubin-Lions theorem give directly
\eqref{co93-d}, so that we can pass to the limit
in the \rhs\ of \eqref{x14ssi}. To take
the limit of the \lhs, on account of the first
\eqref{y12}, it is sufficient to prove that 
\begin{equation} \label{x15}
  \calA(v\ssi) \to \calA(v) \quext{at least weakly in }\,L^2(0,T;V).
\end{equation}  
Since weak convergence surely holds in $\LDH$, it is then
sufficient to prove uniform boundedness in $\LDV$.
With this aim, we compute
\begin{align} \no
  & \nabla \Big( a(v\ssi)\Delta v\ssi + \frac{a'(v\ssi)}2 |\nabla v\ssi|^2 \Big)\\
  \label{x16}  
    & \mbox{}~~~~~
    = a'(v\ssi) \nabla v\ssi \Delta v\ssi + a(v\ssi) \nabla \Delta v\ssi
    + \frac{a''(v\ssi)}2 | \nabla v\ssi |^2 \nabla v\ssi
    + a'(v\ssi) D^2 v\ssi \nabla v\ssi ,
\end{align}  
and, using \eqref{y12}, \eqref{hpa1}, and standard embedding 
properties of Sobolev spaces, it is a standard procedure
to verify that the \rhs\ is uniformly bounded in 
$\LDH$ (and, consequently, so is the left). This concludes 
the proof.
\end{proof}
%



\section{The 6th order problem}
\label{sec:6th}

We start by introducing the concept 
of {\sl weak solution}\ to the sixth order problem
associated with system \eqref{CH1}-\eqref{neum-intro}:
\bede\label{def:weaksol6th}
 Let $\delta>0$ and $\epsi\ge 0$. 
 Let us consider the {\rm 6th order problem} 
 given by the system
 \begin{align}\label{CH1w}
   & u_t + A w = 0, \quext{in }\,V',\\
  \label{CH2w}
   & w = \delta A^2 u + \calA(u) + \xi - \lambda u + \epsi u_t,
    \quext{in }\,V',\\
  \label{CH3w}
   & \xi \in \fzw(u)
 \end{align}
 together with the initial condition
 \begin{equation}\label{init}
   u|_{t=0}=u_0,
   \quext{a.e.~in }\,\Omega.
 \end{equation}
 A (global in time) {\rm weak solution} 
 to the 6th order problem\/
 \eqref{CH1w}-\eqref{init}
 is a triplet $(u,w,\xi)$, with
 \begin{align}\label{regou}
   & u\in \HUVp\cap L^\infty(0,T;W) \cap L^2(0,T;H^3(\Omega)),
    \qquad \epsi u\in \HUH,\\
   \label{regoFu}
   & F(u) \in L^\infty(0,T;L^1(\Omega)),\\
   \label{regofu}
   & \xi \in L^2(0,T;V'),\\
  \label{regow}
   & w\in L^2(0,T;V),
 \end{align}
 satisfying\/ \eqref{CH1w}-\eqref{CH3w} a.e.~in~$(0,T)$
 together with~\eqref{init}. 
\edde
\noindent%
We can then state the main result of this section:
\bete\label{teoesi6th}
 Let us assume\/ \eqref{hpa1}-\eqref{hpf2}. Let
 $\epsi\ge 0$ and $\delta>0$. Moreover, let us suppose
 that
 \begin{equation}\label{hpu0}
   u_0\in W, \quad
    F(u_0)\in L^1(\Omega), \quad
    (u_0)\OO \in (-1,1),
 \end{equation}
 where $(u_0)\OO$ is the spatial mean of $u_0$. 
 Then, the sixth order problem admits one and only one
 weak solution.
 %
 %
 %
\ente
\noindent%
The proof of the theorem will be carried out in several steps,
presented as separate subsequences.
\beos\label{rem:mean}
 We observe that the last condition in \eqref{hpu0},
 which is a common assumption when dealing with Cahn-Hilliard equations
 with constraints (cf.~\cite{KNP} for more details),
 does not simply follow from the requirement $F(u_0)\in L^1(\Omega)$.  
 Indeed, $F$ may be bounded over $[-1,1]$, as it 
 happens, for instance, with the logarithmic 
 potential~\eqref{logpot}. In that case,  
 $F(u_0)\in L^1(\Omega)$ simply means $-1 \le u_0 \le 1$ 
 almost everywhere and, without the last \eqref{hpu0},
 we could have initial data that coincide almost everywhere 
 with either of the pure states $\pm1$. However, solutions that assume 
 (for example) the value $+1$ in a set of strictly positive 
 measure cannot be considered, at least in our regularity
 setting. Indeed, if $|\{u=1\}|>0$, then regularity~\eqref{regofu} 
 (which is crucial for passing to the limit in our
 approximation scheme) is broken, because $f(r)$ is {\sl unbounded}\/ 
 for $r\nearrow +1$ and $\xi$ is nothing else than 
 a relaxed version of $f(u)$. 
\eddos


\subsection{Approximation and local existence}
\label{subsec:appr}

First of all, we introduce a suitably approximated 
statement. The monotone function $f_0$ is regularized
by taking a family $\{f\ssi\}$, $\sigma\in(0,1)$,
defined as in Subsection~\ref{sec:weak}. 
Next, we regularize $u_0$ by singular perturbation,
similarly as before (cf.~\eqref{co93}). Namely,
we take $u\zzs$ as the solution to the elliptic problem
\begin{equation}\label{defiuzzd}
  u\zzs + \sigma A u\zzs = u_0,
\end{equation}
and we clearly have, by Hilbert elliptic
regularity results,
\begin{equation}\label{regouzzd}
  u\zzs \in D(A^2)
   \quad\perogni \sigma\in(0,1).
\end{equation}
Other types of approximations of the initial datum are possible,
of course. The choice \eqref{defiuzzd}, beyond its 
simplicity, has the advantage that it preserves the mean value.

\smallskip

\noindent%
{\bf Approximate problem.}~~For $\sigma\in(0,1)$, 
we consider the problem
\begin{align}\label{CH1appr}
  & u_t + A w = 0,\\  
 \label{CH2appr}
  & w = \delta A^2 u + \calA(u) + f\ssi(u) - \lambda u + (\epsi+\sigma) u_t,\\
 \label{inisd}
  & u|_{t=0}=u\zzs,
   \quext{a.e.~in }\,\Omega.
\end{align}
We shall now show that it admits at least one local 
in time weak solution. Namely, there holds
the following
\bele\label{teo:loc:appr}
 Let us assume\/ \eqref{hpa1}-\eqref{hpf2}.
 Then, for any $\sigma\in(0,1)$,
 there exist $T_0\in(0,T]$ 
 (possibly depending on $\sigma$)
 and a pair $(u,w)$ with
 \begin{align}\label{regovsd}
   & u\in H^1(0,T_0;H)
    \cap L^\infty(0,T_0;W) \cap L^2(0,T_0;D(A^2)),\\
  \label{regowsd}
   & w \in L^2(0,T_0;W),
 \end{align}
 such that \eqref{CH1appr}-\eqref{CH2appr} hold 
 a.e.~in~$(0,T_0)$ and the initial condition~\eqref{inisd}
 is satisfied.
\enle 
\begin{proof}
The theorem will be proved by using the Schauder fixed point theorem.
We take 
\begin{equation}\label{defiBR}
  B_R:=\big \{ v\in L^2(0,T_0;W)\cap L^4(0,T_0;W^{1,4}(\Omega)) : 
   \| v \|_{L^2(0,T_0;W)} + \| v \|_{L^4(0,T_0;W^{1,4}(\Omega))}\le R \big\},
\end{equation}
for $T_0$ and $R$ to be chosen below. 
Then, we take $\baru \in B_R$ and 
consider the problem given by \eqref{inisd} and
\begin{align}\label{CH1schau}
  & u_t + A w = 0, \quext{in }\,H,\\
 \label{CH2schau}
  & w = \delta A^2 u + \calA(\baru) + f\ssi(u) - \lambda u + (\epsi+\sigma) u_t,
   \quext{in }\,H.
\end{align}
Then, 
as $\baru\in B_R$ is fixed, we 
can notice that
\begin{equation}\label{conto21c}
  \| \calA(\baru) \|_{L^2(0,T_0;H)}^2
   \le c \big( \| \baru \|_{L^2(0,T_0;W)}^2
    + \| \baru \|_{L^4(0,T_0;W^{1,4}(\Omega))}^4 \big)
    \le  Q(R).
\end{equation}
Here and below, $Q$ denotes a computable function, 
possibly depending on $\sigma$, defined
for any nonnegative value of its argument(s) and
increasingly monotone in (each of) its argument(s).

As we substitute into \eqref{CH1schau} the expression 
for $w$ given by \eqref{CH2schau}
and apply the inverse operator
$(\Id + (\epsi + \sigma )A )^{-1}$, 
we obtain a parabolic equation in $u$ which is 
linear up to the Lipschitz perturbation $f\ssi(u)$. 
Hence, owing to the
regularity \eqref{conto21c} of the forcing term, 
to the regularity \eqref{regouzzd}
of the initial datum, and to the standard Hilbert
theory of linear parabolic equations, 
there exists a unique pair $(u,w)$ solving the problem 
given by \eqref{CH1schau}-\eqref{CH2schau} and
the initial condition \eqref{inisd}. Such a pair
satisfies the regularity properties 
\eqref{regovsd}-\eqref{regowsd} (as it will also
be apparent from the forthcoming a priori estimates).
We then note as $\calK$ the map
such that $\calK: \baru \mapsto u$. To conclude the proof
we will have to show the following three properties:\\[2mm]
{\sl (i)}~~$\calK$ takes its values in $B_R$;\\[1mm]
{\sl (ii)}~~$\calK$ is continuous with respect to the $L^2(0,T_0;W)$
 and the $L^4(0,T_0;W^{1,4}(\Omega))$ norms;\\[1mm]
{\sl (iii)}~~$\calK$ is a compact map.\\[2mm]
To prove these facts, we perform a couple of a priori estimates.
To start, we test \eqref{CH1schau} by $w$ and 
\eqref{CH2schau} by $u_t$ (energy estimate). This gives
\begin{align} \no
  & \ddt \bigg( \frac\delta2 \| A u \|^2 
   + \io \Big( F\ssi(u) - \frac\lambda2 u^2 \Big) \bigg)
   + (\epsi+\sigma) \| u_t \|^2
   + \| \nabla w \|^2\\
  \label{contox11}
  & \mbox{}~~~~~ 
  = - \big( \calA(\baru) , u_t \big)
   \le \frac\sigma2 \| u_t \|^2
    + \frac1{2\sigma}\| \calA(\baru) \|^2 
\end{align} 
and, after integration in time, 
the latter term can be estimated using
\eqref{conto21c}. Next, we observe
that, thanks to \eqref{propFsigma}, we have
\begin{equation}\label{contox12}
  \frac\delta2 \| A u \|^2 
  + \io \Big( F\ssi(u) - \frac\lambda2 u^2 \Big)
   \ge \eta \| u \|_W^2 - c,   
\end{equation}
for some $\eta>0$, $c\ge 0$ independent of $\sigma$ and 
for all $u$ in $W$.
%
%
%
Thus, \eqref{contox11} provides the bounds
\begin{equation}\label{boundx11}
  \| u \|_{L^\infty(0,T_0;W)}
   + \| u_t \|_{L^2(0,T_0;H)}
   + \| \nabla w \|_{L^2(0,T_0;H)} \le Q\big(R,T_0,\| u\zzs \|_W \big).
\end{equation}
Next, testing \eqref{CH2schau} by $A^2 u$ and 
performing some standard computations (in particular,
the terms $(\calA(\baru),A^2 u)$ and $(f\ssi(u),A^2u)$
are controlled by using \eqref{conto21c}, H\"older's
and Young's inequalities, and the Lipschitz continuity
of $f\ssi$), we obtain the further bound
\begin{equation}\label{st21}
  \| A^2 u \|_{L^2(0,T_0;H)}
   \le Q\big(R,T_0,\|u\zzs\|_{W}\big).
\end{equation}
Hence, estimates \eqref{boundx11} and \eqref{st21}
and a standard application of the Aubin-Lions lemma
permit to see that the range of $\calK$ is 
relatively compact both in $L^2(0,T_0;W)$
and in $L^4(0,T_0;W^{1,4}(\Omega))$.
Thus, {\sl (iii)}\/ follows.

\medskip

Concerning {\sl (i)}, we can now simply observe that,
by \eqref{boundx11},
\begin{equation}\label{st31}
  \| u \|_{L^2(0,T_0;W)}
   \le T_0^{1/2} \| u \|_{L^\infty(0,T_0;W)}
   \le T_0^{1/2} Q\big(R,T_0,\|u\zzs\|_{W}\big).
\end{equation}
whence the \rhs\ can be made smaller than $R$ if $T_0$ is chosen
small enough. A similar estimate works also for the
$L^4(0,T_0;W^{1,4}(\Omega))$-norm since $W\subset W^{1,4}(\Omega)$
continuously. Thus, also {\sl (i)}\ is proved. 

\medskip

Finally, to prove condition {\sl (ii)}, 
we first observe that, 
if $\{\baru_n\}\subset B_R$ 
converges strongly to $\baru$
in $L^2(0,T_0;W)\cap L^4(0,T_0;W^{1,4}(\Omega))$,
then, using proper weak compactness theorems,
it is not difficult to prove that
\begin{equation}\label{conto31}
  \calA(\baru_n)\to \calA(\baru)
   \quext{weakly in }\,L^2(0,T_0;H).
\end{equation}
Consequently, if $u_n$ (respectively $u$) 
is the solution to \eqref{CH1schau}-\eqref{CH2schau}
corresponding to $\baru_n$ (respectively $\baru$), 
then estimates \eqref{boundx11}-\eqref{st21} 
hold for the sequence $\{u_n\}$ with a function $Q$ independent
of $n$. Hence, standard weak compactness arguments
together with the Lipschitz continuity of $f\ssi$ 
permit to prove that
\begin{equation}\label{st33}
  u_n=\calK(\baru_n) \to u=\calK(\baru)
   \quext{strongly in }\,L^2(0,T_0;W) \cap L^4(0,T_0;W^{1,4}(\Omega)),
\end{equation}
i.e., condition {\sl (ii)}. The proof of the lemma
is concluded.
\end{proof}


\subsection{A priori estimates}
\label{sec:apriori}

In this section we will show that the local solutions
constructed in the previous section satisfy uniform
estimates with respect both to the approximation 
parameter $\sigma$ and to 
the time $T_0$. By standard extension methods this 
will yield a global in time solution
(i.e., defined over the whole of $(0,T)$)
in the limit. However, to avoid
technical complications, we will directly assume that the 
approximating solutions are already defined over $(0,T)$. 
Of course, to justify this, we will have to take care
that all the constants appearing in the forthcoming
estimates be independent of $T_0$.
To be precise, in the sequel 
we will note by $c>0$ a computable
positive constant (whose value can vary on occurrence)
independent of all approximation parameters
(in particular of $T_0$ and $\sigma$) and also of 
the parameters $\epsi$ and $\delta$.

\smallskip

\noindent%
{\bf Energy estimate.}~~%
First, integrating \eqref{CH1appr} in space
and recalling \eqref{defiuzzd},
we obtain the {\sl mass conservation}\/ property
\begin{equation}\label{consmedie}
  (u(t))\OO = (u\zzs)\OO
   = (u_0)\OO.
\end{equation}
Next, we can test \eqref{CH1appr} by $w$, \eqref{CH2appr} by $u_t$
and take the difference, arriving at
\begin{equation}\label{conto41}
  \ddt \calE\ssid(u)
   + \| \nabla w \|^2
   + (\epsi+\sigma) \| u_t \|^2
   = 0,
\end{equation}
where the ``approximate energy'' $\calE\ssid(u)$ is defined as
\begin{equation}\label{defiEssid}
  \calE\ssid(u)=\io \Big( \frac\delta2 | A u |^2 
   + \frac{a(u)}2 |\nabla u|^2 + F\ssi(u)
   - \frac{\lambda}2u^2 \Big).
\end{equation}
Actually, it is clear that the high regularity of
approximate solutions (cf.~\eqref{regovsd}-\eqref{regowsd})
allows the integration by parts necessary to write \eqref{conto41}
(at least) almost everywhere in time. 
Indeed, all single terms in \eqref{CH2appr} lie in 
$\LDH$ and the same holds for the test function $u_t$.

Then, we integrate \eqref{conto41} in time and notice
that, by \eqref{propFsigma},
\begin{equation}\label{Essicoerc}
  \calE\ssid(u) \ge \eta \big(
   \delta \| u \|_W^2 + \| u \|_V^2 
   \big) - c   \quad\perogni t\in(0,T).
\end{equation}
Consequently, \eqref{conto41} provides the bounds
\begin{align} \label{st41}
  & \| u \|_{\LIV} + \delta^{1/2} \| u \|_{\LIW} + (\epsi+\sigma)^{1/2} 
   \| u_t \|_{\LDH} \le c,\\
 \label{st43}
  & \| \nabla w \|_{\LDH} \le c,\\
 \label{st44}
  & \| F\ssi(u) \|_{L^\infty(0,T;L^1(\Omega))} \le c,
\end{align}  
where it is worth stressing once more
that the above constants $c$ 
neither depend explicitly on $\delta$ nor on $\epsi$.

\smallskip

\noindent%
{\bf Second estimate.}~~%
We test \eqref{CH2appr} by $u-u\OO$, $u\OO$ denoting
the (constant in time) spatial mean of $u$. Integrating
by parts the term $\calA(u)$, we obtain
\begin{align}\no
  & \delta \| A u \|^2 
   + \io a(u) | \nabla u |^2
   + \io f\ssi(u)\big( u - u\OO \big)\\
 \label{conto51}
  & \mbox{}~~~~~ 
   \le \big( w + \lambda u - (\epsi+\sigma) u_t, u - u\OO \big)
   - \io \frac{a'(u)}2 | \nabla u |^2 ( u - u\OO )
\end{align}
and we have to estimate some terms. First
of all, we observe that there exists
a constant $c$, depending on the (assigned
once $u_0$ is fixed) value of $u\OO$, but
{\sl independent of $\sigma$},
such that
\begin{equation}\label{conto52}
  \io f\ssi(u)\big( u - u\OO \big)
   \ge \frac12 \| f\ssi(u) \|_{L^1(\Omega)} - c.
\end{equation}
To prove this inequality, one basically
uses the monotonicity of $f\ssi$ and the fact
that $f\ssi(0)=0$ (cf.~\cite[Appendix]{MZ}
or \cite[Third a priori estimate]{GMS}
for the details). Next, by 
\eqref{hpa2}, the function 
$r\mapsto a'(r)(r-u\OO)$ is uniformly bounded,
whence
\begin{equation}\label{conto53}
  - \io \frac{a'(u)}2 | \nabla u |^2 ( u - u\OO )
   \le c \| \nabla u \|^2.
\end{equation}
Finally, using that $(w\OO+\lambda u\OO, u-u\OO)=0$ 
since $w\OO+\lambda u\OO$ is constant with respect to space variables,
and applying the Poincar\'e-Wirtinger inequality,
\begin{align}\no
  & \big( w + \lambda u - (\epsi+\sigma) u_t, u - u\OO \big)
   = \big( w - w\OO + \lambda (u-u\OO) - (\epsi+\sigma) u_t, u - u\OO \big)\\
 \no
  & \mbox{}~~~~~
   \le c \| \nabla w \| \| \nabla u \|
    + c \| \nabla u \|^2
    + c (\epsi + \sigma) \| u_t \| \| \nabla u \|  \\
 \label{conto54}
  & \mbox{}~~~~~   
   \le c\big( \| \nabla w \| 
   + (\epsi + \sigma) \| u_t \|
   + 1 \big),
\end{align}
the latter inequality following from 
estimate~\eqref{st41}.

Thus, squaring \eqref{conto51}, using
\eqref{conto52}-\eqref{conto54}, and integrating
in time, we arrive after recalling 
\eqref{st41}, \eqref{st43} at
\begin{equation} \label{st51}
  \| f\ssi(u) \|_{L^2(0,T;L^1(\Omega))} \le c.
\end{equation}
Next, integrating \eqref{CH2appr} with respect 
to space variables (and, in particular, integrating
by parts the term $\calA(u)$), using
\eqref{st51}, and recalling \eqref{st43}, we obtain
(still for $c$ independent of $\epsi$ and $\delta$)
\begin{equation} \label{st52}
  \| w \|_{L^2(0,T;V)} \le c.
\end{equation}

\noindent%
{\bf Third estimate.}~~%
We test \eqref{CH2appr} by $Au$.
Using the monotonicity of 
$f\ssi$ and \eqref{hpa1}, it is
not difficult to arrive at
\begin{equation}\label{conto61}
  \frac{\epsi+\sigma}2\ddt \| \nabla u \|^2
   + \delta \| \nabla A u \|^2 
   + \frac{\agiu}2 \| A u \|^2
  \le \big( \nabla w + \lambda \nabla u , \nabla u \big)
  + c \| \nabla u \|_{L^4(\Omega)}^4.
\end{equation}
Using the continuous embedding 
$H^{3/4}(\Omega)\subset L^4(\Omega)$
(so that, in particular, 
$H^{7/4}(\Omega)\subset W^{1,4}(\Omega)$)
together with the interpolation inequality 
\begin{equation}\label{new-interp}
  \| v \|_{H^{3/4}(\Omega)} 
   \le \| v \|^{3/8}_{H^3(\Omega)}  \| v \|^{5/8}_{H^1(\Omega)} 
   \quad \perogni v \in H^{3/4}(\Omega),
\end{equation}
and recalling estimate \eqref{st41},
the last term is treated as follows:
\begin{equation}\label{conto62}
  c \| \nabla u \|_{L^4(\Omega)}^4 
   \le c \| u \|_{H^3(\Omega)}^{3/2}
    \| u \|_{V}^{5/2}
   \le \frac\delta2 \| \nabla A u \|^2
    + c(\delta).
\end{equation}
Note that the latter constant $c(\delta)$
is expected to explode as $\delta\searrow 0$ but, on
the other hand, is independent of $\sigma$.
Next, noting that
\begin{equation}\label{conto63}
  \big( \nabla w + \lambda \nabla u , \nabla u \big)
   \le c \big( \| \nabla u \|^2 + \| \nabla w \|^2 ),
\end{equation}
from \eqref{conto61} we readily deduce
\begin{equation} \label{st61}
  \| u \|_{L^2(0,T;H^3(\Omega))} 
   \le c(\delta).
\end{equation}
A similar (and even simpler)
argument permits to check that it is also
\begin{equation} \label{st62}
  \| \calA(u) \|_{L^2(0,T;H)} \le c(\delta).
\end{equation}
Thus, using \eqref{st41}, \eqref{st52}, 
\eqref{st61}-\eqref{st62} and comparing terms in
\eqref{CH2appr}, we arrive at
\begin{equation} \label{st63}
  \| f\ssi(u) \|_{L^2(0,T;V')} \le c(\delta).
\end{equation}
%


\subsection{Limit $\boldsymbol \sigma\searrow 0$}
\label{sec:sigma}

We now use the machinery introduced in 
Subsection~\ref{sec:weak} to take the 
limit $\sigma\searrow 0$ in 
\eqref{CH1appr}-\eqref{CH2appr}. For convenience,
we then rename as $(u\ssi,w\ssi)$ the solution. Then, 
recalling estimates \eqref{st41}-\eqref{st44},
\eqref{st52} and \eqref{st61}-\eqref{st63},
and using the Aubin-Lions compactness lemma,
we deduce
\begin{align} \label{conv41}
  & u\ssi \to u \quext{strongly in }\, 
   C^0([0,T];H^{2-\epsilon}(\Omega)) \cap 
   L^2(0,T;H^{3-\epsilon}(\Omega)),\\
 \label{conv42}
  & u\ssi \to u \quext{weakly star in }\, H^1(0,T;V') \cap 
   L^\infty(0,T;W) \cap 
   L^2(0,T;H^3(\Omega)),\\
 \label{conv42b}
  & (\epsi+\sigma) u_{\sigma,t} \to \epsi u_t 
   \quext{weakly in }\, L^2(0,T;H),\\
 \label{conv43}
  & w\ssi \to w \quext{weakly in }\, \LDV,\\
 \label{conv44}
  & f\ssi(u\ssi) \to \xi \quext{weakly in }\, \LDVp,
\end{align}  
for suitable limit functions $u,w,\xi$, where $\epsilon>0$ is 
arbitrarily small. It is readily checked that the above 
relations (\eqref{conv41} in particular) are strong
enough to guarantee that
\begin{equation} \label{conv45}
  \calA(u\ssi) \to \calA(u), \quext{strongly in }\, \LDH.
\end{equation}
This allows us to take the limit $\sigma\searrow 0$
in \eqref{CH1appr}-\eqref{inisd} (rewritten for $u\ssi,w\ssi$) 
and get 
\begin{align}\label{CH1delta}
  & u_t + A w = 0, \quext{in }\,V',\\
 \label{CH2delta}
  & w = \delta A^2 u + \calA(u) + \xi - \lambda u + \epsi u_t,
   \quext{in }\,V',\\
 \label{iniz-delta}
  & u|_{t=0} = u_0
   \quext{a.e.~in }\,\Omega.
\end{align}
To identify $\xi$, we observe that,
thanks to \eqref{conv41}, \eqref{conv44}, and 
Lemma~\ref{limimono} applied with the choices of 
$X=V$, $X'=V'$, $\calB\ssi=f\ssi$, $\calB=\fzw$,
$v\ssi=u\ssi$, $v=u$ and $\xi\ssi=f\ssi(u\ssi)$,
it follows that 
\begin{equation} \label{incldelta}
  \xi\in \fzw(u).
\end{equation}
Namely, $\xi$ is identified with respect to the 
weak (duality) expression of the function $f_0$.
This concludes the proof of Theorem~\ref{teoesi6th}
for what concerns existence.


\subsection{Uniqueness}
\label{sec:uniq6th}

To conclude the proof of Theorem~\ref{teoesi6th}, it remains to 
prove uniqueness. To this purpose, we write both \eqref{CH1w} and
\eqref{CH2w} for a couple of solutions $(u_1,w_1,\xi_1)$, $(u_2,w_2,\xi_2)$, 
and take the difference. This gives
\begin{align}\label{CH1d0}
  & u_t + A w = 0, \quext{in }\,V',\\
 \no
  & w = \delta A^2 u - a(u_1) \Delta u
   - \big( a(u_1) - a(u_2) \big) \Delta u_2
   - \frac{a'(u_1)}2 \big( | \nabla u_1 |^2 - | \nabla u_2 |^2 \big)\\
 \label{CH2d0}
  & \mbox{}~~~~~~~~~~
   - \frac{a'(u_1) - a'(u_2)}2 | \nabla u_2 |^2 
   + \xi_1 - \xi_2 - \lambda u + \epsi u_t,
   \quext{in }\,V',
\end{align}
where we have set $(u,w,\xi):=(u_1,w_1,\xi_1)-(u_2,w_2,\xi_2)$.
Then, we test \eqref{CH1d0} by $A^{-1}u$, \eqref{CH2d0} by $u$,
and take the difference. Notice that, indeed, $u$ has zero mean value
by \eqref{consmedie}. Thus, the operator $A^{-1}$ makes 
sense since $A$ is bijective from $W_0$ to $H_0$, 
the subscript $0$ indicating the zero-mean condition.

A straighforward computation involving use of standard
embedding properties of Sobolev spaces then gives
\begin{align}\no
  & \Big(- a(u_1) \Delta u - \big( a(u_1) - a(u_2) \big) \Delta u_2
   - \frac{a'(u_1)}2 \big( | \nabla u_1 |^2 - | \nabla u_2 |^2 \big)
   - \frac{a'(u_1) - a'(u_2)}2 | \nabla u_2 |^2 , u \Big) \\
 \label{uniq22} 
  & \mbox{}~~~~~
   \le Q\big( \| u_1 \|_{L^\infty(0,T;W)},\| u_2 \|_{L^\infty(0,T;W)} \big)
   \| u \|_W \| u \| 
\end{align}
and we notice that the norms inside the function $Q$ are controlled
thanks to \eqref{regou}. Thus, also on account of the monotonicity of 
$\fzw$, we arrive at
\begin{align}\no
  & \ddt\Big( \frac12 \| u \|_{V'}^2 + \frac\epsi2 \| u \|^2 \Big) 
   + \delta \| A u \|^2
   \le c \| u \|_W \| u \| + \lambda \| u \|^2\\
 \label{uniq23}
  & \mbox{}~~~~~
   \le c \| u \|_W^{4/3} \| u \|_{V'}^{2/3} 
   \le \frac\delta2 \| A u \|^2 + c(\delta) \| u \|_{V'}^2,  
\end{align}
where, to deduce the last two inequalities, we used
the interpolation inequality $\| u \| \le \| u \|_{V'}^{2/3} 
\| u \|_W^{1/3}$ (note that $(V,H,V')$ form 
a Hilbert triplet, cf., e.g., \cite[Chap.~5]{BrAF})
together with Young's inequality
and the fact that the function
$\| \cdot \|_{V'} + \| A \cdot \|$
is an equivalent norm on $W$. Thus, the thesis of 
Theorem~\ref{teoesi6th} follows by applying Gronwall's 
lemma to \eqref{uniq23}.


\section{From the 6th order to the 4th order model}
\label{sec:6thto4th}

In this section, we analyze the behavior of solutions to
the 6th order problem as $\delta$ tends to $0$. To start with,
we specify the concept of weak solution in the 4th order case:
\bede\label{def:weaksol4th}
 Let $\delta=0$ and $\epsi\ge 0$. 
 Let us consider the\/ {\rm 4th order problem} given by
 the system
 \begin{align}\label{CH1w4th}
   & u_t + A w = 0, \quext{in }\,V',\\
  \label{CH2th}
   & w = \calA(u) + f(u) + \epsi u_t,
    \quext{in }\,H,
 \end{align}
 together with the initial condition~\eqref{init}.
 A\/ (global in time) {\rm weak solution} 
 to the 4th order problem \eqref{CH1w4th}-\eqref{CH2th}, 
 \eqref{init} is a pair $(u,w)$, with
 \begin{align}\label{regou4}
   & u\in \HUVp\cap L^\infty(0,T;V) \cap L^2(0,T;W),
    \qquad \epsi u\in \HUH,\\
   \label{regoFu4}
   & F(u) \in L^\infty(0,T;L^1(\Omega)),\\
   \label{regofu4}
   & f_0(u) \in L^2(0,T;H),\\
  \label{regow4}
   & w\in L^2(0,T;V),
 \end{align}
 satisfying\/ \eqref{CH1w4th}-\eqref{CH2th} a.e.~in~$(0,T)$
 together with\/ \eqref{init}.
\edde
\noindent%
\bete\label{teo6thto4th}
 Let us assume\/ \eqref{hpa1}-\eqref{hpf2} together with
 \begin{equation}\label{aconcave}
   a \quext{is concave on }\,[-1,1].
 \end{equation}
 Let also $\epsi\ge 0$ and let, for all $\delta\in(0,1)$, 
 $u\zzd$ be an initial datum satisfying\/ \eqref{hpu0}.
 Moreover, let us suppose
 \begin{equation}\label{convuzzd}
   u\zzd\to u_0 \quext{strongly in }\,V, \qquad
    \calE\dd(u\zzd)\to \calE_0(u_0),
    \quext{where }\,(u_0)\OO\in(-1,1).
 \end{equation}
 Let, for any $\delta\in (0,1)$, $(u\dd,w\dd,\xi\dd)$ be a 
 weak solution to the 6th order system 
 in the sense of\/
 {\rm Definition~\ref{def:weaksol6th}}. Then, we have that, 
 up to a (nonrelabelled) subsequence of $\delta\searrow 0$,
 \begin{align}\label{co4th11}
   & u\dd \to u \quext{weakly star in }\,\HUVp \cap \LIV \cap \LDW,\\
   \label{co4th12}
   & \epsi u\dd \to \epsi u \quext{weakly in }\,\HUH,\\
   \label{co4th13}
   & w\dd \to w \quext{weakly in }\,\LDV,\\
   \label{co4th13b}
   & \delta u\dd \to 0 \quext{strongly in }\,L^2(0,T;H^3(\Omega)),\\
   \label{co4th14}
   & \xi\dd \to f_0(u) \quext{weakly in }\,\LDVp,    
 \end{align}
 and $(u,w)$ is a weak solution to the 4th order problem.
\ente
\noindent%
\begin{proof}
The first part of the proof consists in repeating the 
``Energy estimate'' and the ``Second estimate''
of the previous section. In fact, we could avoid this procedure
since we already noted that the constants appearing
in those estimates were independent of $\delta$. 
However, we choose to perform once more the estimates
working directly on the 6th order problem
(rather than on its approximation) for various
reasons. First, this will show that the estimates
do not depend on the chosen regularization scheme. 
Second, the procedure has an independent interest
since we will see that the use of ``weak'' 
subdifferential operators still permits to rely
on suitable integration by parts formulas and
on monotonicity methods. Of course, many passages,
which were trivial in the ``strong'' setting,
need now a precise justification. 
Finally, in this way we are able to prove, as a 
byproduct, that any solution to the 6th order system
satisfies an energy {\sl equality}\/ (and not just
an inequality). Actually, this property may be useful for
addressing the long-time behavior of the system.

\smallskip

\noindent%
{\bf Energy estimate.~~}%
As before, we would like to test \eqref{CH1w} by $w\dd$,
\eqref{CH2w} by $u_{\delta,t}$, and take the difference. To justify
this procedure, we start observing that $w\dd\in L^2(0,T;V)$ 
by \eqref{regow}. Actually, since \eqref{CH1w} is in 
fact a relation in $L^2(0,T;V')$, the use of $w\dd$ as 
a test function makes sense. The problem, instead, arises when
working on \eqref{CH2w} and, to justify the estimate, we 
can just consider the (more difficult) case $\epsi=0$. 

%
%

Then, it is easy to check that the assumptions of 
Lemma~\ref{lemma:ipp} are satisfied. In particular, we have 
\eqref{x12} thanks to \eqref{regou}.
Hence, \eqref{x14} gives 
\begin{equation}\label{en-11}
  \duavb{u_{\delta,t},\calA(u\dd)} 
   = \frac12 \ddt \io a(u\dd)|\nabla u\dd|^2, 
    \quext{a.e.~in }\,(0,T).
\end{equation}
Thus, it remains to show that
\begin{equation}\label{en-12}
  \duavg{u_{\delta,t},\delta A^2 u\dd + \xi\dd}
   = \ddt \io \Big( \frac\delta2 |A u\dd|^2
    + F(u\dd) \Big),
    \quext{a.e.~in }\,(0,T).
\end{equation}
To prove this, we observe that 
\begin{equation}\label{comparis}
  \delta A^2 u\dd + \xi\dd \in \LDV \quad
   \perogni \delta\in(0,1).
\end{equation}
Actually, we already noted above that $w\dd$, $\calA(u\dd)$ 
lie in $\LDV$. Since $u\dd\in \LDV$ by \eqref{regou}
and we assumed $\epsi=0$, \eqref{comparis} simply
follows by comparing terms in \eqref{CH2w}.
Thus, the duality on the \lhs\ of \eqref{en-12} makes sense. 
Moreover, as we set
\begin{equation}\label{en-13}
  \calJ\dd(v):= \io \Big( \frac\delta2 |A v|^2 + F(v) \Big),
\end{equation}
then a direct computation permits to check that
\begin{equation}\label{en-14}
  \delta A^2 u\dd + \xi\dd
   \in \de \calJ\dd ( u\dd )
   \quext{a.e.~in }\,(0,T).
\end{equation}
Indeed, by definition
of $H$-subdifferential, this corresponds to 
the relation
\begin{equation}\label{en-15}
  \duavb { \delta A^2 u\dd + \xi\dd , v - u\dd }
   \le \calJ\dd ( v ) - \calJ\dd ( u\dd )
  \quad \perogni v\in H,  
\end{equation}
and it is sufficient to check it for $v\in V$ since for
$v\in V\setminus H$ the \rhs\ is $+\infty$ and consequently
the relation is trivial. However, for $v\in V$, 
\eqref{en-15} follows by definition of the relaxed
operator $\fzw$. Thanks to \eqref{en-14},
\eqref{en-12} is a then a direct consequence of 
inequality \eqref{ipepardiff} of Lemma~\ref{BResteso}.

Thus, the above procedure permits to see that
(any) weak solution $(u\dd,w\dd,\xi\dd)$ to the
6th order problem satisfies the energy {\sl equality}
\begin{equation}\label{energy-6th}
  \ddt \calE\dd(u(t)) 
   + \| \nabla w(t) \|^2
   + \epsi \| u_t(t) \|^2 = 0
\end{equation}
for almost all $t\in[0,T]$. As a consequence, we 
get back the first two convergence relations 
in \eqref{co4th11} as
well as \eqref{co4th12}. Moreover, we have
\begin{equation}\label{6to4-01}
  \| \nabla w\dd \|_{\LDH} \le c.
\end{equation}

\smallskip

\noindent%
{\bf Second estimate.~~}%
Next, to get \eqref{co4th13} and \eqref{co4th14},
we essentially need to repeat the ``Second estimate''
of the previous section. Indeed, we see that $u\dd-(u\dd)\OO$ is
an admissible test function in \eqref{CH1w}. However,
we now have to obtain an estimate of $\xi\dd$ from
the duality product
\begin{equation}\label{6to4-21}
  \duavb{\xi\dd, u\dd - (u\dd)\OO}.
\end{equation}
Actually, if $\xi\dd=\xi\dda+\xi\dds$ is the Lebesgue
decomposition of the {\sl measure}\ $\xi\dd$
given in Theorem~\ref{teobrezis},
then, noting that for all $t\in [0,T]$ we have  
$u\dd(t)\in W\subset C^0(\barO)$, we can write
\begin{equation}\label{6to4-22}
  \duavb{\xi\dd(t), u\dd(t) - (u\dd)\OO}
   = \io \xi\dda(t) \big(u\dd(t) - (u\dd)\OO\big)\,\dix
   + \io \big(u\dd(t) - (u\dd)\OO \big) \dixi_{\delta,s}(t).
\end{equation}
Next, we notice that, as a direct
consequence of assumption~\eqref{convuzzd},
\begin{equation}\label{unifsep}
  \esiste\mu\in(0,1):~~
   -1+\mu \le (u\zzd)\OO \le 1-\mu,
   \quad\perogni \delta\in(0,1),
\end{equation}
where $\mu$ is independent of $\delta$.
In other words, the spatial means $(u\zzd)\OO$ are
uniformly separated from $\pm1$. 
Then, recalling \eqref{bre2} and proceeding 
as in \eqref{conto52}, we have
\begin{equation}\label{conto52dd}
  \io \xi\dda(t) \big(u\dd(t) - (u\dd)\OO\big)\,\dix
   \ge \frac12\| \xi\dda(t) \|_{L^1(\Omega)} - c, 
\end{equation}
where $c$ does not depend on $\delta$. 

On the other hand, let us note as
$\deriv\!\xi\dds=\phi\dds \dixis$ the {\sl polar}\
decomposition of $\xi\dds$, where $|\xi\dds|$ is 
the {\sl total variation}\/ of $\xi\dds$
(cf., e.g., \cite[Chap.~6]{Ru}).
Then, introducing the bounded linear functional 
$\calS\dd:C^0(\barO)\to \RR$ given by
\begin{equation}\label{deficalS}
  \calS\dd(z):= \ibaro z \,\dixi_{\delta,s}
\end{equation}
using, e.g., \cite[Thm.~6.19]{Ru}, 
and recalling \eqref{bre3},
we can estimate the norm of $\calS$ as 
follows:
\begin{align}\no
  |\xi\dds|(\barO)
  & = \ibaro \dixis 
   = \| \calS\dd \|_{\calM(\barO)}\\
 \no 
  & = \sup\Big\{\ibaro z\, \dixi_{\delta,s},~
     z\in C^0(\overline\Omega),~
     z(\barO)\subset[-1,1] \Big\}\\
 \label{6to4-23}
  & = \duav{\xi\dd,u\dd} - \io \xi\dda u\dd
   = \ibaro u\dd \,\dixi_{\delta,s}
   = \ibaro u\dd\phi\dds\,\dixis,
\end{align}
where we also used that $u\dd\in C^0(\barO)$.
Comparing terms, it then follows
\begin{equation}\label{6to4-24}
  u\dd = \phi\dds,
   \quad |\xi\dds|-\text{a.e.~in $\barO$}.
\end{equation}
Then, since is clear that
\begin{equation}\label{6to4-25}
  u\dd=\pm 1 \implica \frac{u\dd-(u\dd)\OO}{|u\dd-(u\dd)\OO|}=\pm1,
\end{equation}
coming back to \eqref{6to4-23} we deduce
\begin{equation}\label{6to4-26}
  \ibaro \dixis
   = \ibaro \phi\dds\frac{u\dd-(u\dd)\OO}{|u\dd-(u\dd)\OO|}\,\dixis
   \le c \ibaro \phi\dds \big( u\dd-(u\dd)\OO \big)\,\dixis
   = c \ibaro \big( u\dd-(u\dd)\OO \big) \, \dixi\dds.
\end{equation}
Here we used again in an essential way the 
uniform separation property \eqref{unifsep}.

Collecting \eqref{6to4-22}-\eqref{6to4-26}, 
we then have
\begin{equation}\label{6to4-21b}
  \duavb{\xi\dd, u\dd - (u\dd)\OO}
   \ge \frac12 \| \xi\dda(t) \|_{L^1(\Omega)} 
   + \eta \ibaro \dixis - c,
\end{equation}
for some $c\ge0$, $\eta>0$ independent of $\delta$.
On the other hand, mimicking \eqref{conto51}-\eqref{conto54},
we obtain
\begin{equation}\label{6to4-27}
  \delta \| A u\dd \|^2
   + \agiu \| \nabla u\dd \|^2
   + \duavb{\xi\dd, u\dd - (u\dd)\OO}
  \le c \big(\| \nabla w\dd \|
  + \epsi \| u_{\delta,t} \| + 1 \big),
\end{equation}
whence squaring, integrating in time,
and using \eqref{co4th12} and \eqref{6to4-01},
we obtain that the function
\begin{equation}\label{6to4-28}
 t \mapsto \| \xi\dda(t) \|_{L^1(\Omega)} 
  + \io \dixist
  \quext{is bounded in $L^2(0,T)$, independently of $\delta$.}
\end{equation}
Integrating now \eqref{CH2w} in space, we deduce
\begin{equation}\label{6to4-29}
  \io w\dd 
   = \frac12 \io a'(u\dd) | \nabla u\dd |^2
    + \io \xi\dd
    - \lambda (u\dd)\OO,
\end{equation}
whence
\begin{equation}\label{6to4-29b}
 \bigg| \io w\dd \bigg|
  \le c \bigg( \| \nabla u\dd \|^2
    + \| \xi\dda(t) \|_{L^1(\Omega)} 
  + \ibaro \dixist + 1 \bigg).
\end{equation}
Thus, squaring, integrating in time, and recalling 
\eqref{6to4-01} and \eqref{6to4-28}, we finally
obtain \eqref{co4th13}.

\smallskip

\noindent%
{\bf Key estimate.}~~%
To take the limit $\delta>0$, we have to provide a 
bound on $\calA(u\dd)$ independent of $\delta$. 
This will be obtained by means of the following
integration by parts formula due to Dal Passo, 
Garcke and Gr\"un (\cite[Lemma 2.3]{DpGG}):
\bele\label{lemma:dpgg}
 Let $h\in W^{2,\infty}(\RR)$ and $z\in W$. Then,
 \begin{align} \no
   & \io h'(z) |\nabla z|^2 \Delta z
    = -\frac13 \io h''(z) |\nabla z|^4\\
  \label{byparts}
   & \mbox{}~~~~~
    + \frac23 \io h(z) \big( |D^2 z|^2 - | \Delta z|^2 \big)
    + \frac23 \iga h(z) II( \nabla z ),
 \end{align}  
 where $II(\cdot)$ denotes the second fundamental
 form of $\Gamma$. 
\enle
\noindent%
We then test \eqref{CH2w} by $Au\dd$ in the duality between
$V'$ and $V$. This gives the relation
\begin{equation} \label{conto71}  
  \frac\epsi2\ddt \| \nabla u\dd \|^2
   + \delta \| \nabla A u\dd \|^2
   +\big( \calA(u\dd), A u\dd \big)
   + \duavg{\xi\dd,Au\dd}
  = \big(\nabla w\dd, \nabla u\dd\big)
    + \lambda \| \nabla u\dd \|^2
\end{equation}
and some terms have to be estimated. First, we note 
that
\begin{equation} \label{conto71b}  
  \big( \calA(u\dd), Au\dd \big)
   = \Big( a(u\dd)\Delta u\dd + \frac{a'(u\dd)}2 |\nabla u\dd|^2,
    \Delta u\dd \Big).
\end{equation}
Thus, using Lemma~\ref{lemma:dpgg} with the choice 
of $h(\cdot)=a'(\cdot)/2$, we obtain
\begin{align} \no
  & \big( \calA(u\dd), Au\dd \big)
  = \io a(u\dd) | \Delta u\dd |^2
   + \frac13 \io a(u\dd) \big( |D^2 u\dd|^2 - | \Delta u\dd|^2 \big)\\
 \label{conto45} 
  & \mbox{}~~~~~~~~~~ 
  - \frac16 \io a''(u\dd)|\nabla u\dd|^4
     + \frac13 \iga a(u\dd) II(\nabla u\dd).
\end{align}
Let us now point out that, being $\Gamma$ smooth, we can estimate
\begin{equation} \label{conto46}
  \frac13 \bigg| \iga a(u\dd)  II(\nabla u\dd) \bigg|
   \le c \| \nabla u\dd \|_{L^2(\Gamma)}^2
   \le \omega \| A u\dd \|_{W}^2
   + c_\omega \| u\dd \|^2,
\end{equation}  
for small $\omega>0$ to be chosen below,
the last inequality following from the continuity of 
the trace operator (applied to $\nabla u$)
from $H^s(\Omega)$ into $L^2(\Gamma)$ for $s\in(1/2,1)$
and the compactness of the embedding 
$W\subset H^{1+s}(\Omega)$ for $s$ in the same range.

Thus, using the {\sl concavity}\/ assumption \eqref{aconcave}
on $a$ and the fact that $|u_\delta|\le 1$ almost everywhere
in $(0,T)\times\Omega)$,
we get
\begin{equation}\label{elliptic3}
  \big( \calA(u\dd), Au\dd \big)
  \ge \eta \| A u\dd \|^2
   - c,
\end{equation}
for proper strictly positive constants $\eta$ and $c$,
both independent of $\delta$.
Next, we observe that, by \eqref{incldelta} and 
Lemma~\ref{BSesteso}, we obtain $\duavg{\xi\dd,Au\dd}\ge 0$.
Finally, we have
\begin{equation}\label{conto73}
  - (\nabla w\dd, \nabla u\dd)
   \le c \| \nabla w\dd \| \| \nabla u\dd \|,
\end{equation}
and the \rhs\ is readily estimated thanks 
to \eqref{co4th12} and \eqref{co4th13}.

Thus, on account of \eqref{elliptic3}, 
integrating \eqref{conto71} in time,
we readily obtain 
the last of \eqref{co4th11} as
well as \eqref{co4th13b}. Moreover, since
$-1\le u\dd\le 1$ almost everywhere, we have
for free
\begin{equation} \label{st-key}
  \| u\dd \|_{L^\infty((0,T)\times\Omega)}\le 1.
\end{equation}
Thus, using the Gagliardo-Nirenberg inequality 
\eqref{ineq:gn}, we have also
\begin{equation} \label{conv54}
  u\dd \to u \quext{weakly in }\, 
   L^4(0,T;W^{1,4}(\Omega)).
\end{equation}  
This readily entails
\begin{equation} \label{conv55}
  \calA(u\dd) \to \calA(u) \quext{weakly in }\, 
   L^2(0,T;H).
\end{equation}  
Thus, a comparison of terms in \eqref{CH2w} gives also
\begin{equation} \label{conv56}
  \xi\dd \to \xi \quext{weakly in }\, 
   L^2(0,T;V').
\end{equation}  
Then, we can take the limit $\delta\searrow 0$ in
\eqref{CH1w} and get \eqref{CH1w4th}. On the other
hand, if we take the limit of \eqref{CH2w}, we
obtain
\begin{equation} \label{CH2provv}
   w = \calA(u) + \xi - \lambda u + \epsi u_t
\end{equation}  
and we have to identify $\xi$. Actually, \eqref{conv56},
the strong convergence $u\dd\to u$ in $\LDV$ 
(following from \eqref{co4th11} and the Aubin-Lions lemma)
and Lemma~\ref{limimono} permit to show that
\begin{equation} \label{incldelta2}
  \xi\in \fzw(u) \quext{a.e.~in }\,(0,T).
\end{equation}  
On the other hand, a comparison argument in
\eqref{CH2provv} permits to see that 
$\xi\in \LDH$, whence, thanks to \eqref{betavsbetaw2},
we obtain that $\xi(t)=f_0(u(t))\in H$ for a.e.~$t\in(0,T)$.
This concludes the proof of Theorem~\ref{teo6thto4th}.
\end{proof}


\section{Analysis of the fourth order problem}
\label{sec:4th}

In this section, we will prove existence of
a weak solution to Problem~\eqref{CH1}-\eqref{neum-intro}
in the fourth order case $\delta =0$ by means of a direct
approach not relying on the 6th order approximation. 
This will allow us to consider a general function 
$a$ (without the concavity assumption \eqref{aconcave}).
More precisely, we have the following
\bete\label{teo:4th}
 Let assumptions\/ \eqref{hpa1}-\eqref{hpf2} hold,
 let $\epsi\ge 0$ and let
 \begin{equation}\label{hpu0-4}
   u_0\in V, \quad
    F(u_0)\in L^1(\Omega), \quad
    (u_0)\OO \in (-1,1).
 \end{equation}
 Then, there exists\/ {\rm at least} one weak
 solution to the 4th order problem, in the sense
 of\/ {\rm Definition~\ref{def:weaksol4th}.}
 %
 %
 %
 %
 %
\ente
\noindent%
The rest of the section is devoted to the proof of the 
above result, which is divided into several steps.

\smallskip

\noindent%
{\bf Phase-field approximation.}~~For 
$\sigma\in(0,1)$, we consider the system
\begin{align}\label{CH1-4ap}
  & u_t + \sigma w_t + A w = 0,\\  
 \label{CH2-4ap}
  & w = \calA(u) + f\ssi(u) - \lambda u + (\epsi+\sigma) u_t.
\end{align}
This will be endowed with the initial conditions
\begin{equation}\label{init-4ap}
  u|_{t=0} = u\zzs, \qquad
   w|_{t=0} = 0.
\end{equation}
Similarly as before (compare with \eqref{defiuzzd}),
we have set
\begin{equation}\label{defiuzzs}
  u\zzs + \sigma A^2 u\zzs = u_0
\end{equation}
and, by standard elliptic regularity, we have that
\begin{equation}\label{propuzzs}
  u\zzs \in H^5(\Omega)\subset C^{3+\alpha}(\barO) 
   \quext{for }\,\alpha\in(0,1/2),
   \qquad \dn u\zzs=\dn A u\zzs=0,~~\text{on }\,\Gamma.
\end{equation}
Moreover, of course, $u\zzs\to u_0$ in a suitable sense
as $\sigma\searrow0$.

\smallskip

\noindent%
{\bf Fixed point argument.}~~We now prove existence of 
a local solution to the phase-field
approximation by a further Schauder fixed point argument. 
Namely, we introduce the system
\begin{align}\label{CH1-4pf}
  & u_t + \sigma w_t + A w = 0,\\  
 \label{CH2-4pf}
  & \barw = -a(\baru) \Delta u - \frac{a'(\baru)}2 | \nabla \baru |^2 
   + f\ssi(\baru) - \lambda \baru + (\epsi + \sigma) u_t, \qquad
  \dn u = 0~~\text{on }\,\Gamma,
\end{align}
which we still endow with the condition \eqref{init-4ap}.
Here, $f\ssi$ is chosen as in \eqref{defifsigma}.

Next, we set 
\begin{equation}\label{deficalU}
  \calU:=\left\{ u\in C^{0,1+\alpha} ([0,T_0]\times \barO):~
   u|_{t=0}=u\zzs,~  
   \| u \|_{C^{0,1+\alpha}}\le 2R \right\},
\end{equation}
where $R:=\max\{1,\|u\zzs\|_{C^{1+\alpha}(\barO)}\}$
and $T_0$ will be chosen at the end of the argument.
It is clear that $R$ depends in fact on $\sigma$
(so that the same will happen for $T_0$). This dependence 
is however not emphasized here.
For the definition of the parabolic H\"older spaces 
used in this proof we refer the reader to 
\cite[Chap.~5]{Lu}, whose notation is adopted.
Moreover, in the sequel,
in place of $C^{0,\alpha}([0,T_0]\times \barO)$
(and similar spaces) we will just write $C^{0,\alpha}$,
for brevity. We then also define 
\begin{equation}\label{deficalW}
  \calW:=\left\{ w\in C^{0,\alpha}:~
   w|_{t=0}=0,~  
   \| w \|_{C^{0,\alpha}}\le R \right\},
\end{equation}
where $R$ is, for simplicity, the same number 
as in \eqref{deficalU}.

Then, choosing $(\baru,\barw)$ in $\calU\times\calW$
and inserting it in \eqref{CH2-4pf}, we observe that,
by the Lipschitz regularity of $a$ (cf.~\eqref{hpa1})
and standard multiplication properties of
H\"older spaces, there exists a computable
monotone function $Q$, also depending on $\sigma$, 
but independent of the time $T_0$, such that
\begin{equation}\label{prop-norme}
  \|a(\baru)\|_{C^{0,\alpha}}
  + \big\|a'(\baru)|\nabla \baru|^2\big\|_{C^{0,\alpha}}
  + \|f\ssi(\baru)\|_{C^{0,\alpha}}
  \le Q(R).
\end{equation}
Thanks to \cite[Thm.~5.1.21]{Lu}, then there exists one and only one
solution $u$ to \eqref{CH2-4pf} with the first initial condition
\eqref{init-4ap}. This solution satisfies
\begin{equation}\label{regou-4pf}
  \| u \|_{C^{1,2+\alpha}}
   \le Q(R).
\end{equation}
Substituting then $u_t$ in \eqref{CH1-4pf} and applying the same
theorem of \cite{Lu} to this equation with the second initial condition
\eqref{init-4ap}, we then obtain one and only one solution
$w$, with
\begin{equation}\label{regow-4pf}
  \| w \|_{C^{1,2+\alpha}}
   \le Q(R).
\end{equation}
We then note as $\calT$ the map
such that $\calT: (\baru,\barw) \mapsto (u,w)$. As before,
we need to show that:\\[2mm]
{\sl (i)}~~$\calT$ takes its values in $\calU\times\calW$;\\[1mm]
{\sl (ii)}~~$\calT$ is continuous with respect to the 
 $C^{0,1+\alpha}\times C^{0,\alpha}$ norm
 of $\calU\times\calW$;\\[1mm]
{\sl (iii)}~~$\calT$ is a compact map.\\[2mm]
First of all let us prove {\sl (i)}. We just refer to the component
$u$, the argument for $w$ being analogous and in fact simpler. 
We start observing that, if $u\in \Pi_1 (\calT (\calU\times \calW))$,
$\Pi_1$ denoting the projection on the first component, then
\begin{equation}\label{i-11}
  \| u(t) \|_{C^\alpha(\barO)}
   \le \| u_0 \|_{C^\alpha(\barO)}
   + \int_0^t \| u_t(s) \|_{C^\alpha(\barO)} \,\dis 
   \le R + T_0 Q(R), 
   \quad \perogni t\in[0,T_0],
\end{equation}
which is smaller than $2R$ if $T_0$ is chosen suitably.

Next, using the continuous embedding (cf.~\cite[Lemma~5.1.1]{Lu})
\begin{equation} \label{contemb}
  C^{1,2+\alpha} \subset 
   C^{1/2}([0,T_0];C^{1+\alpha}(\barO))
  \cap C^{\alpha/2}([0,T_0];C^2(\barO)),
\end{equation}
we obtain that, analogously,
\begin{equation}\label{i-12}
  \| \nabla u(t) \|_{C^\alpha(\barO)}
   \le \| \nabla u_0 \|_{C^\alpha(\barO)}
   + T_0^{1/2} \| u \|_{C^{1/2}([0,T_0];C^{1+\alpha}(\barO))}
   \le R + T_0^{1/2} Q(R).
\end{equation}
Hence, passing to the supremum for $t\in[0,T_0]$, we see
that the norm of $u$ in $C^{0,1+\alpha}$ can be made
smaller than $2R$ if $T_0$ is small enough. 
Thus, {\sl (i)}\ is proved.

\medskip

Let us now come to {\sl (iii)}. As before, we just deal with the 
component $u$. Namely, on account of \eqref{regou-4pf},
we have to show that the space $C^{1,2+\alpha}$ is compactly
embedded into $C^{0,1+\alpha}$. Actually, by
\eqref{contemb} and using standard compact inclusion 
properties of H\"older spaces, this relation is proved easily.
Hence, we have {\sl (iii)}.

\medskip

Finally, we have to prove {\sl (ii)}. This property
is however straighforward. Actually, taking 
$(\baru_n,\barw_n)\to (\baru,\barw)$ in 
$\calU\times \calW$, we have that the corresponding
solutions $(u_n,w_n)=\calT(\baru_n,\barw_n)$ 
are bounded in the sense of \eqref{regou-4pf}-\eqref{regow-4pf}
uniformly in $n$. Consequently, a standard weak 
compactness argument, together with
the uniqueness property for the initial value problems
associated to \eqref{CH1-4pf} and to \eqref{CH2-4pf},
permit to see that the {\sl whole sequence}\/
$(u_n,w_n)$ converges to a unique limit point 
$(u,w)$ solving \eqref{CH1-4pf}-\eqref{CH2-4pf}
with respect to the limit data $(\baru,\barw)$. Moreover, by the 
compactness property proved in {\sl (iii)}, this
convergence holds with respect to the original topology
of $\calU\times \calW$. This proves that
$(u,w) = \calT (\baru,\barw)$, i.e., {\sl (ii)}\
holds. 

\medskip

\noindent%
{\bf A priori estimates.}~~For any $\sigma>0$, we have obtained
a local (i.e., with a final time $T_0$ depending on $\sigma$) 
solution to \eqref{CH1-4ap}-\eqref{CH2-4ap} with the 
initial conditions \eqref{init-4ap}. To emphasize the
$\sigma$-dependence, we will note it by $(u\ssi,w\ssi)$ 
in the sequel. To let $\sigma\searrow 0$, we now
devise some of a priori estimates
uniform both with respect to $\sigma$ and with respect to
$T_0$. As before, 
this will give a global solution in the limit and, to
avoid technicalities, we can
directly work on the time interval $[0,T]$. 
Notice that the high regularity of $(u\ssi,w\ssi)$ 
gives sense to all the calculations performed below
(in particular, to all the integrations by parts).
That said, we repeat the ``Energy estimate'', exactly as
in the previous sections. This now gives
\begin{align} \label{st11ap}
  & \| u\ssi \|_{\LIV} + \| F\ssi(u\ssi) \|_{L^\infty(0,T;L^1(\Omega))} \le c,\\
 \label{st12ap}
  & (\sigma+\epsi)^{1/2} \| u_{\sigma,t} \|_{\LDH} \le c,\\
 \label{st13ap}
  & \sigma^{1/2} \| w\ssi \|_{L^\infty(0,T;H)} 
    + \| \nabla w\ssi \|_{\LDH} \le c.
\end{align}  
Next, working as in the ``Second estimate'' of
Subsection~\ref{sec:apriori}, we obtain the analogue 
of \eqref{st51} and \eqref{st52}.

To estimate $f\ssi(u\ssi)$ in $H$, we 
now test \eqref{CH2-4ap} by $f\ssi(u\ssi)$, to get
\begin{align}\no
  & \frac{\epsi+\sigma}2 \ddt \io F\ssi(u\ssi)
   + \io \Big( a(u\ssi) f\ssi'(u\ssi) 
      + \frac{a'(u\ssi)}2 f\ssi(u\ssi) \Big) | \nabla u\ssi |^2       
   + \| f\ssi(u\ssi) \|^2 \\ 
 \label{conto51-4th}
  & \mbox{}~~~~~ 
   = \big( w\ssi + \lambda u\ssi, f\ssi(u\ssi) \big),
\end{align}
and it is a standard matter to estimate the \rhs\ by using
the last term on the \lhs, H\"older's and Young's inequalities, 
and properties \eqref{st11ap} and \eqref{st52}. 
Now, we notice that, thanks to \eqref{goodmono},
\begin{equation}\label{4th-21}
  a(r) f\ssi'(r) + \frac{a'(r)}2 f\ssi(r)
   \ge \agiu f\ssi'(r) - c | f\ssi(r) |
   \ge \frac{\agiu}2 f\ssi'(r) - c
   \quad \perogni r\in [-2,2],
\end{equation}
with the last $c$ being independent of $\sigma$.
On the other hand, for $r\not\in [-2,2]$
we have that $a'(r)=0$ by \eqref{hpa2}.
Hence, also thanks to \eqref{st11ap}, the second term 
on the \lhs\ of \eqref{conto51-4th} can be controlled. 
We then arrive at
\begin{equation} \label{st21ap}
  \| f\ssi(u\ssi) \|_{L^2(0,T;H)} \le c.
\end{equation}
The key point is represented by the next estimate, which is
used to control the second space derivatives of $u$. To do 
this, we have to operate a change of variable, first.
Namely, we set
\begin{equation} \label{defiphi}
  \phi(s):=\int_0^s a^{1/2} (r) \, \dir, \qquad
   z\ssi:=\phi(u\ssi)
\end{equation}
and notice that, by \eqref{hpa1}-\eqref{hpa2}, $\phi$ is monotone
and Lipschitz together with its inverse. Then, by \eqref{st11ap},
\begin{equation} \label{st21ap2}
  \| z\ssi \|_{L^\infty(0,T;V)} \le c
\end{equation}
and it is straighforward to realize that \eqref{CH2-4ap}
can be rewritten as
\begin{equation} \label{CH2-z}
  w\ssi = - \phi'(u\ssi) \Delta z\ssi 
     + f\ssi\circ\phi^{-1}(z\ssi) - \lambda u\ssi 
     + (\epsi + \sigma) u_{\sigma,t}, 
   \qquad \dn u\ssi = 0~~\text{on }\,\Gamma.
\end{equation}
By the H\"older continuity of $u\ssi$ up to its second space
derivatives and the Lipschitz continuity of $a$ and $a'$
(cf.~\eqref{hpa1}-\eqref{hpa2}), $-\Delta z\ssi$ is also H\"older
continuous in space.
Thus, we can use it as a test function in \eqref{CH2-z}.
Using the monotonicity of $f\ssi$ and $\phi^{-1}$,
and recalling \eqref{st21ap}, we then easily obtain
\begin{equation} \label{st31ap}
  \| z\ssi \|_{L^2(0,T;W)} \le c.
\end{equation}

\smallskip

\noindent%
{\bf Passage to the limit.}~~As a consequence 
of \eqref{st11ap}-\eqref{st13ap}, 
\eqref{st51}-\eqref{st52} and
\eqref{st21ap}, we have
\begin{align} \label{co11ap}
  & u\ssi \to u \quext{weakly star in }\,
   \HUVp \cap \LIV,\\
 \label{co12ap} 
  & (\sigma+\epsi) u_{\sigma,t} \to \epsi u_t 
   \quext{weakly in }\, \LDH,\\
 \label{co13ap} 
  & f\ssi(u\ssi) \to \barf 
   \quext{weakly in }\, \LDH,\\
 \label{co14ap}
  & w\ssi \to w
   \quext{weakly in }\, \LDV,\\
 \label{co15ap}
  & u_{\sigma,t} + \sigma w_{\sigma,t} 
  \to u_t \quext{weakly in }\, \LDVp,
\end{align}  
for suitable limit functions $u$, $w$ and $\barf$.
Here and below, all convergence relations have to 
be intended to hold up to (nonrelabelled) subsequences
of $\sigma\searrow0$. Now, by the Aubin-Lions lemma,
we have
\begin{equation} \label{co21ap}
  u\ssi \to u \quext{strongly in }\, \CZH
   \quext{and a.e.~in }\,Q.
\end{equation}
Then, \eqref{co13ap} and a standard monotonicity
argument (cf.~\cite[Prop.~1.1]{barbu}) 
imply that $\barf=f(u)$ a.e.~in $Q$.
Furthermore, by \eqref{hpa1}-\eqref{hpa2}
and the generalized Lebesgue's theorem,
we have 
\begin{equation} \label{co21ap2}
  a(u\ssi) \to a(u),~~a'(u\ssi) \to a'(u),~~
   \quext{strongly in }\,
   L^q(Q)~~\text{for all }\,q\in[1,+\infty).
\end{equation}
Analogously, recalling \eqref{st21ap2},
$z\ssi=\phi(u\ssi)\to \phi(u)=:z$,
strongly in $L^q(Q)$ for all $q\in[1,6)$.
Actually, the latter relation holds also weakly 
in $\LDW$ thanks to the bound \eqref{st31ap}. 
Moreover, by \eqref{st21ap2}, \eqref{st31ap} 
and interpolation, we obtain
\begin{equation} \label{co22ap}
  \| \nabla z\ssi \|_{L^{10/3}(Q)} \le c,
\end{equation}
whence, clearly, it is also
\begin{equation} \label{co22ap2}
  \| \nabla u\ssi \|_{L^{10/3}(Q)} \le c.
\end{equation}
As a consequence, being 
\begin{equation} \label{co22ap3}
  - \Delta u\ssi = - \frac1{a^{1/2}(u\ssi)} \Delta z\ssi
   + \frac{a'(u\ssi)}{2a(u\ssi)} | \nabla u\ssi |^2,
\end{equation}
we also have that
\begin{equation} \label{co22ap4}
  \Delta u\ssi \to \Delta u 
   \quext{weakly in }\, L^{5/3}(Q).
\end{equation}
Combining this with \eqref{co11ap} and using the 
generalized Aubin-Lions lemma (cf., e.g., \cite{Si}),
we then arrive at
\begin{equation} \label{co22ap5}
  u\ssi \to u 
   \quext{strongly in }\, L^{5/3}(0,T;W^{2-\epsilon,5/3}(\Omega))
    \cap C^0([0,T];H^{1-\epsilon}(\Omega)),
    \quad \perogni \epsilon > 0,
\end{equation}
whence, by standard interpolation and embedding properties
of Sobolev spaces, we obtain
\begin{equation} \label{co22ap6}
  \nabla u\ssi \to \nabla u 
   \quext{strongly in }\, L^q(Q) \quext{for some }\,q>2.
\end{equation}
Consequently, recalling \eqref{co21ap2},
\begin{equation} \label{co22ap7}
  a'(u\ssi) |\nabla u\ssi|^2 \to a'(u) |\nabla u|^2,
   \quext{say, weakly in }\, L^1(Q).
\end{equation}
This is sufficient to take the limit $\sigma\searrow 0$
in \eqref{CH2-4ap} and get back \eqref{CH2th}.
To conclude the proof, it only remains to show the
regularity \eqref{regou4} for what concerns the second
space derivatives of $u$. Actually, by 
\eqref{st31ap} and the Gagliardo-Nirenberg inequality
\eqref{ineq:gn},
\begin{equation} \label{co22ap8}
  z \in L^2(0,T;W) \cap L^\infty(Q) \subset
   L^4(0,T;W^{1,4}(\Omega)).
\end{equation}
Thus, we have also $u \in L^4(0,T;W^{1,4}(\Omega))$ and, 
consequently, a comparison of terms in \eqref{CH2th}
permits to see that $\Delta u \in L^2(0,T;H)$, whence 
\eqref{regou4} follows from elliptic regularity.
The proof of Theorem~\ref{teo:4th} is concluded.


\section{Further properties of weak solutions}
\label{sec:uniq}


\subsection{Uniqueness for the 4th order problem}
\label{subsec:uniq}

We will now prove that, if the interfacial (i.e., gradient) 
part of the free energy $\calE\dd$ satisfies a 
{\sl convexity}\/ condition (in the viscous case $\epsi>0$)
or, respectively, a {\sl strict convexity}\/ condition 
(in the non-viscous case $\epsi=0$), then the solution 
is unique also in the 4th order case. 
Actually, the stronger assumption
(corresponding to $\kappa>0$ in the statement below)
required in the non-viscous case is needed for the purpose 
of controlling the nonmonotone part of $f(u)$, while in the 
viscous case we can use the term $\epsi u_t$ 
for that aim.

It is worth noting that, 
also from a merely thermodynamical point of view, 
the convexity condition is a rather natural requirement. 
Indeed, it corresponds to asking the second differential 
of $\calE\dd$ to be positive definite, to ensure
that the stationary solutions are dynamically stable 
(cf., e.g., \cite{Su} for more details). 
\bete\label{teouniq}
 Let the assumptions of\ {\rm Theorem~\ref{teo:4th}} hold
 and assume that, in addition,
 \begin{equation} \label{1aconc}
   a''(r)\ge 0, \quad
   \Big(\frac1a\Big)''(r)\le -\kappa,
    \quad\perogni r\in[-1,1],
 \end{equation}  
 where $\kappa>0$ if $\epsi=0$ and $\kappa\ge 0$ 
 if $\epsi> 0$. Then, the 4th order problem admits a
 unique weak solution.
\ente 
\begin{proof} 
Let us denote by $J$ the gradient part of the energy,
i.e.,
\begin{equation} \label{defiJ}
  J:V\to [0,+\infty), \qquad
   J(u):=\io \frac{a(u)}2 | \nabla u |^2.
\end{equation}  
Then, we clearly have 
\begin{equation} \label{Jprime}
  \duavg{J'(u),v} =
    \io \Big( a(u) \nabla u \cdot \nabla v
    + \frac{a'(u)}2 |\nabla u|^2 v \Big).   
\end{equation}  
and we can correspondingly compute the second derivative
of $J$ as
\begin{equation}\label{Jsecond}
  \duavg{J''(u)v,z} = 
  \io \Big( \frac{a''(u)|\nabla u|^2 vz}2
   + a'(u) v \nabla u\cdot\nabla z 
   + a'(u) z \nabla u\cdot\nabla v 
   + a(u) \nabla v\cdot\nabla z \Big).
\end{equation}  
To be more precise, we have that
$J'(u)\in V'$ and $J''(u)\in \calL(V,V')$
at least for $u\in W$ (this may instead not be true
if we only have $u\in V$, due to the quadratic terms
in the gradient). This is however the case for the 
4th order system since for any weak solution
we have that $u(t) \in W$ at least for a.e.~$t\in(0,T)$.

From \eqref{Jsecond}, we then have in particular
\begin{align} \no
  \duavg{J''(u)v,v} 
  & = \io \Big( \frac{a''(u)|\nabla u|^2 v^2}2
    + 2 a'(u) v \nabla u\cdot\nabla v 
    + a(u) | \nabla v |^2 \Big)\\
 \label{Jvv}    
  & \ge \io \Big( a(u) - \frac{2 a'(u)^2}{a''(u)} \Big)
   | \nabla v|^2, 
\end{align}  
whence the functional $J$ is convex, at least when restricted
to functions $u$ such that
\begin{equation} \label{doveJconv}
  u \in W, \quad u(\Omega)\subset [-1,1],
\end{equation}  
provided that $a$ satisfies
\begin{equation} \label{aaprimo}
  a(r)a''(r) - 2a'(r)^2 \ge 0
   \quad\perogni r\in [-1,1].
\end{equation}  
Noting that
\begin{equation} \label{1asecondo}
  \Big(\frac1a\Big)''
   = \frac{2(a')^2-aa''}{a^3},
\end{equation}  
we have that $J$ is (strictly) convex 
if $1/a$ is (strictly) concave, i.e.,
\eqref{1aconc} holds (cf.~also \cite[Sec.~3]{DNS}
for related results). Note that, 
in deducing the last inequality in \eqref{Jvv},
we worked as if it was $a''>0$. However, if
$a''(r)=0$ for some $r$, then also $a'(r)$
has to be $0$ due to \eqref{aaprimo}. So, this
means that in the set $\{u=r\}$ the first two
summands in the \rhs\ of the first line of 
\eqref{Jvv} identically vanish.

\smallskip

That said, let us write both \eqref{CH1w} and
\eqref{CH2w} for a couple of solutions $(u_1,w_1)$, $(u_2,w_2)$, 
and take the difference.
Setting $(u,w):=(u_1,w_1)-(u_2,w_2)$, we obtain
\begin{align}\label{CH1d}
  & u_t + A w = 0,\\
 \label{CH2d}
  & w = J'(u_1) - J'(u_2) + f(u_1) - f(u_2) + \epsi u_t.
\end{align}
%
%
Then, we can test \eqref{CH1d} by
$A^{-1}u$, \eqref{CH2d} by $u$, and take the difference.
Indeed, $u=u_1-u_2$ has zero mean value
by \eqref{consmedie}. We obtain
\begin{equation} \label{contod1}
  \frac12 \ddt \Big( \| u \|_{V'}^2 + \epsi \| u \|^2 \Big)
   + \duavg{J'(u_1) - J'(u_2),u}
   + \big( f(u_1) - f(u_2), u \big) = 0   
\end{equation}  
and, using the convexity of $J$ 
coming from \eqref{1aconc} and the 
$\lambda$-monotonicity of $f$
(see~\eqref{hpf1}), we have,
for some function $\xi$ belonging to $W$ 
a.e.~in time and taking its values in $[-1,1]$,
\begin{equation} \label{contod2}
  \frac12 \ddt \Big( \| u \|_{V'}^2 + \epsi \| u \|^2 \Big)
   + \kappa \| \nabla u \|^2
  \le \frac12 \ddt \Big( \| u \|_{V'}^2 + \epsi \| u \|^2 \Big)
   + \duavg{J''(\xi) u, u}
   \le \lambda \| u \|^2.
\end{equation}   
Thus, in the case $\epsi>0$ (where it may be $\kappa=0$),
we can just use Gronwall's Lemma. Instead, if 
$\epsi=0$ (so that we assumed $\kappa>0$),
by the Poincar\'e-Wirtinger inequality we
have
\begin{equation} \label{contod3}
  \lambda \| u \|^2  
  \le \frac\kappa2 \| \nabla u \|^2
   + c \| u \|_{V'}^2,
\end{equation}
and the thesis follows again by applying 
Gronwall's lemma to \eqref{contod2}.
\end{proof}


\subsection{Additional regularity}
\label{sec:add}

We prove here parabolic regularization properties
of the solutions to the 4th order system
holding in the case of a convex energy functional.
An analogous result would hold also for the 6th order 
system under general conditions on $a$ since the
bilaplacean in that case dominates the lower order 
terms (we omit the details).
\bete\label{teoreg}
 Let the assumptions of\/ {\rm Theorem~\ref{teouniq}} hold.
 Then, the solution satisfies the additional
 regularity property
 \begin{equation} \label{add-reg}
   \| u \|_{L^\infty(\tau,T;W)} 
    + \| u \|_{L^\infty(\tau,T;W^{1,4}(\Omega))}
    \le Q(\tau^{-1}) \quad \perogni \tau>0,
 \end{equation}  
 where $Q$ is a computable monotone function whose expression
 depends on the data of the problem and, in particular,
 on $u_0$.
\ente
\begin{proof}
The proof is based on a further a priori estimate, which
has unfortunately a formal character in the present 
regularity setting. To justify it, one should proceed
by regularization. For instance, a natural choice would be
that of refining the fixed point argument leading to 
existence of a weak solution (cf.~Sec.~\ref{sec:4th}) 
by showing (e.g., using a bootstrap regularity argument)
that, at least locally in time, the solution lies 
in higher order H\"older spaces. We leave the details
to the reader.

That said, we test \eqref{CH1w} by $w_t$ and subtract
the result from the time derivative of \eqref{CH2w}
tested by $u_t$. We obtain
\begin{equation} \label{contoe1}
  \frac 12 \ddt \| \nabla w \|^2
  + \frac\epsi2 \ddt \| u_t \|^2
   + \duavg{J''(u) u_t, u_t}
   + \io f'(u) u_t^2 \le 0.
\end{equation}  
Then, by convexity of $J$,
\begin{equation} \label{contoe2}
  \duavg{J''(u) u_t, u_t}
   \ge \kappa \| u_t \|_{\LDV}^2.
\end{equation}  
On the other hand, the
$\lambda$-monotonicity of $f$ gives
\begin{equation} \label{contoe3}
  \io f'(u) u_t^2 
   \ge - \lambda \| u_t \|_{\LDH}^2
\end{equation}  
and, if $\epsi=0$ (so that $\kappa>0$),
we have as before
\begin{equation} \label{contoe3-b}
  - \lambda \| u_t \|_{\LDH}^2
  \ge - \frac\kappa2 \| u_t \|_{\LDV}^2
   - c \| u_t \|_{\LDVp}^2.
\end{equation}  
Thus, recalling the first of \eqref{regou}
and applying the {\sl uniform}\/ Gronwall lemma
(cf.~\cite[Lemma~I.1.1]{Te}),
it is not difficult to infer
\begin{equation} \label{ste1}
  \| \nabla w \|_{L^\infty(\tau,T;H)} 
   + \epsi^{1/2} \| u_t \|_{L^\infty(\tau,T;H)} 
   + \kappa \| u_t \|_{L^2(\tau,T;V)} 
   \le Q(\tau^{-1}) \quad \perogni \tau>0.
\end{equation}  
Next, testing \eqref{CH2w} by $u-u\OO$ and proceeding
as in the ``Second estimate'' of Subsection~\ref{sec:apriori},
but taking now the essential supremum as time
varies in $[\tau,T]$, we arrive at
\begin{equation} \label{ste2}
  \| w \|_{L^\infty(\tau,T;V)} 
   + \| f(u) \|_{L^\infty(\tau,T;L^1(\Omega))} 
   \le Q(\tau^{-1}) \quad \perogni \tau>0.
\end{equation}  
Thus, thanks to \eqref{ste2}, we can test \eqref{CH2th} by
$-\Delta z$, with $z=\phi(u)$ (cf.~\eqref{defiphi}). Proceeding
similarly with Section~\ref{sec:4th} (but taking now the 
supremum over $[\tau,T]$ rather than integrating in time),
we easily get \eqref{add-reg}, which concludes the proof.
\end{proof}


\subsection{Energy equality}
\label{sec:long}

As noted in Section~\ref{sec:6thto4th},
any weak solution to the 6th order system satisfies 
the energy {\sl equality} \eqref{energy-6th}.
We will now see that the same property holds 
also in the {\sl viscous}\/ 4th order
case (i.e., if $\delta=0$ and $\epsi>0$). 
More precisely, we can prove the 
\bepr\label{prop:energy}
 Let the assumptions of\/ {\rm Theorem~\ref{teo:4th}}
 hold and let $\epsi>0$. Then, any weak solution to the
 4th order system satisfies the\/ {\rm integrated} energy
 equality
 \begin{equation}\label{energy-4th-i}
   \calE_0(u(t)) = \calE_0(u_0) 
    - \int_0^t \big( \| \nabla w(s) \|^2
    - \epsi \| u_t(s) \|^2 \big) \dis 
    \quad \perogni t\in[0,T].
 \end{equation}
\empr
\begin{proof}
As before, we proceed by testing \eqref{CH1w4th} by $w$,
\eqref{CH2th} by $u_t$ and taking the difference. As
$u_t\in \LDH$ and $f_0(u)\in \LDH$
(cf.~\eqref{regou4} and \eqref{regofu4}),
then the integration by parts
\begin{equation} \label{co91}
  \big(f(u),u_t\big)
   = \ddt \io F(u), \quext{a.e.~in }\,(0,T)
\end{equation}  
is straighforward (it follows directly
from \cite[Lemma~3.3, p.~73]{Br}). 

Moreover, in view of \eqref{regou4},
assumption \eqref{x11} of Lemma~\ref{lemma:ipp}
is satisfied. Hence, by \eqref{x14},
we deduce that
\begin{equation} \label{co92}
  \int_0^t \big(\calA(u(s)),u_t(s)\big)\,\dis
   = \io \frac{a(u(t))}2 |\nabla u(t)|^2
   - \io \frac{a(u_0)}2 |\nabla u_0|^2.
\end{equation}  
Combining \eqref{co91} and \eqref{co92}, we immediately get
the assert.
\end{proof}
\noindent%
It is worth noting that the energy equality obtained 
above has a key relevance in the investigation 
of the long-time behavior of
the system. In particular, given $m\in(-1,1)$ (the spatial mean
of the initial datum, which is a conserved quantity due
to~\eqref{consmedie}), we can define the {\sl phase space}
\begin{equation} \label{defiXd}
  \calX\ddm:=\big\{u\in V:~\delta u\in W,~F(u)\in L^1(\Omega),~
   u\OO = m\big\}
\end{equation}  
and view the system (both for $\delta>0$ and for $\delta=0$)
as a (generalized) dynamical process in $\calX\ddm$.
Then, \eqref{energy-4th-i} (or its 6th order analogue)
stands at the basis of the so-called {\sl energy method}\/ 
(cf.~\cite{Ba1,MRW}) for proving existence of the 
{\sl global attractor} with respect to the ``strong''
topology of the phase space.
This issue will be analyzed in a forthcoming work.
\beos\label{nonviscous}
 Whether the equality \eqref{energy-4th-i} still holds
 in the nonviscous case $\epsi=0$ seems to be
 a nontrivial question.
 The answer would be positive in case one could prove 
 the integration by parts formula
 \begin{equation} \label{co101}
   \itt\duavb{u_t,\calA(u)+f(u)}
    = \io\Big(\frac{a(u(t))}2 |\nabla u(t)|^2 + F(u(t))\Big)
    - \io\Big(\frac{a(u_0)}2|\nabla u(t)|^2 + F(u_0)\Big),
 \end{equation}  
 under the conditions
 \begin{equation} \label{co102}
   u \in \HUVp \cap L^2(0,T;W) \cap L^\infty(Q), 
    \qquad \calA(u)+f(u) \in \LDV,
 \end{equation}  
 which are satisfied by our solution (in particular
 the latter \eqref{co102} follows by a comparison of 
 terms in \eqref{CH2th}, where now $\epsi=0$).
 Actually, if \eqref{co102} holds, then both hands
 sides of \eqref{co101} make sense. However, devising
 an approximation argument suitable for proving
 \eqref{co101} could be a rather delicate problem.
\eddos



\vspace{15mm}

\noindent%
{\bf First author's address:}\\[1mm]
Giulio Schimperna\\
Dipartimento di Matematica, Universit\`a degli Studi di Pavia\\
Via Ferrata, 1,~~I-27100 Pavia,~~Italy\\
E-mail:~~{\tt giusch04@unipv.it}

\vspace{4mm}

\noindent%
{\bf Second author's address:}\\[1mm]
Irena Paw\l ow\\
Systems Research Institute,\\
Polish Academy of Sciences\\
\mbox{}~~and Institute of Mathematics and Cryptology,\\
Cybernetics Faculty,\\
Military University of Technology,\\
S.~Kaliskiego 2,~~00-908 Warsaw,~~Poland\\ 
E-mail:~~{\tt Irena.Pawlow@ibspan.waw.pl}

\end{document}